\documentclass{article}

\usepackage{amsmath}
\usepackage{latexsym}
\usepackage{graphicx}
\usepackage{psfrag}
\usepackage{amssymb}
\usepackage{amscd}
\usepackage{fullpage}

\def\<{\langle}
\def\>{\rangle}
\def\0{{{\bf 0}}}

\def\CG{{\mathcal G}}
\def\OO{{\mathcal O}}

\def\CA{{\mathcal A}}

\def\CM{{\mathcal M}}

\def\CV{{\mathcal V}}

\def\AA{{\mathbb A}}
\def\CC{{\mathbb C}}

\def\FF{{\mathbb F}}

\def\QQ{{\mathbb Q}}
\def\RR{{\mathbb R}}
\def\TT{{\mathbb T}}
\def\ZZ{{\mathbb Z}}

\def\pp{{\mathfrak p}}
\def\ll{{\mathfrak l}}

\def\ti{{\tilde i}}
\def\tm{{\tilde m}}

\def\tX{{\tilde X}}

\def\tU{{\tilde U}}

\def\tPsi{{\tilde \Psi}}

\newcommand{\spec}{\operatorname{Spec}}

\newcommand{\im}{\operatorname{im}}

\newcommand{\Hom}{\operatorname{Hom}}

\newcommand{\End}{\operatorname{End}}
\newcommand{\BC}{\operatorname{BC}}

\newcommand{\JL}{\operatorname{JL}}
\newcommand{\Lie}{\operatorname{Lie}}
\newcommand{\Frob}{\operatorname{Frob}}

\newcommand{\gal}{\operatorname{Gal}}
\newcommand{\Aut}{\operatorname{Aut}}
\newcommand{\DR}{\operatorname{DR}}
\newcommand{\GL}{\operatorname{GL}}
\newcommand{\GU}{\operatorname{GU}}

\newcommand{\cris}{\mbox{\rm \tiny cris}}
\newcommand{\univ}{\mbox{\rm \tiny univ}}
\newcommand{\et}{\mbox{\rm \tiny \'et}}
\newcommand{\conn}{\mbox{\rm \tiny conn}}
\newcommand{\tor}{\mbox{\rm \tiny tor}}
\newcommand{\sep}{\mbox{\rm \tiny sep}}

	{\begin{list}
		{\noindent\makebox[0mm][r]{\arabic{enumi}.}}
		{\usecounter{enumi} \topsep=1.5mm \itemsep=0mm}
	}
	{\end{list}}

	{\begin{list}
		{\noindent\makebox[0mm][r]{\arabic{enumi}.}}
		{\usecounter{enumi} \topsep=1.5mm \itemsep=-1mm}
	}
	{\end{list}}

\begin{document}%%%%%%%%%%%%%%%%%%%%%%%%%%%%%%%%%%%%%%%%%%%%%%%%%%%%%%%%%
%%%%%%%%%%%%%%%%%%%%%%%%%%%%%%%%%%%%%%%%%%%%%%%%%%%%%%%%%%%%%%%%%%%%%%%%%

% Theorem environments with italic font
%\theoremstyle{plain}
  \newtheorem{theorem}{Theorem}[section]
  \newtheorem{definition}[theorem]{Definition}
  \newtheorem{lemma}[theorem]{Lemma}
  \newtheorem{corollary}[theorem]{Corollary}
  \newtheorem{proposition}[theorem]{Proposition}
  \newtheorem{conjecture}[theorem]{Conjecture}
  \newtheorem{question}[theorem]{Question}
  \newtheorem{problem}[theorem]{Problem}

% Theorem environments with roman font
  \newtheorem{example}[theorem]{Example}
  \newtheorem{remark}[theorem]{Remark}

\newenvironment{proof}{{{\noindent{\it Proof.\ }}}}{{\hfill $\Box$}}

\title{Mazur's principle for U(2,1) Shimura varieties}

\author{David Helm
\footnote{Harvard University; dhelm@math.harvard.edu}}

\maketitle

%\begin{abstract}
Mazur's principle gives a criterion under which an irreducible mod $l$
Galois representation arising from a classical modular form of
level $Np$ (with $p$ prime to $N$) also arises from a classical modular
form of level $N$.  We consider the analogous question for Galois
representations arising from certain unitary Shimura varieties.
In particular, we prove an analogue of Mazur's principle for $U(2,1)$
Shimura varieties.  We also give a conjectural criterion for a Galois
representation arising in the cohomology of a unitary Shimura variety
whose level subgroup is parahoric at $p$ to also arise in the
cohomology of a Shimura variety with ``less level structure at $p$''.

%\end{abstract}

%%%%%%%%%%%%%%%%%%%%%%%%%%%%%%%%%%%%%%%%%%%%%%%%%%%%%%%%%%%%%%%%%%%%%%%%%
\section{Introduction}

A key result in the theory of classical modular forms was the determination
of the optimal level from which a given modular mod $l$ representation could
arise.  A crucial step in this direction was given by Mazur's principle,
which states:

\begin{theorem}[Mazur's Principle (\cite{Ri1}, Theorem 6.1)]
Let $N$ be an integer, let $p$ be a prime not dividing $N$, and let
$\overline{\rho}: \gal(\overline{\QQ}/\QQ) \rightarrow \GL_2(\overline{\FF_l})$
be a Galois representation arising from a modular form on $\Gamma_0(Np)$.
Suppose that:
\begin{itemize}
\item $p \neq l$,
\item $\overline{\rho}$ is absolutely irreducible and unramified at $p$, and
\item $l$ does not divide $p-1$.
\end{itemize}
Then $\overline{\rho}$ arises from a modular form for $\Gamma_0(N)$.
\end{theorem}

Mazur's principle was an important ingredient in
Ribet's proof of his level-lowering theorem~\cite{Ri1}.
More recently, Jarvis~\cite{Ja} gave a generalization of Mazur's principle
to Galois representations attached to Shimura curves over totally real
fields; Rajaei~\cite{Ra} has made use of this generalization to extend 
Ribet's level-lowering theorem to this setting as well.

The analogous problem of determining the optimal levels for mod $l$ 
representations attached to unitary groups remains wide open, and one can
ask if there is an analogue to Mazur's principle in this setting.  Where
Mazur's principle starts with a representation arising from a form
with $\Gamma_0(p)$ level structure at $p$ and attempts to ``remove $p$'' 
from this level structure, such an analogue should start with a 
representation $\overline{\rho}$ arising
in the cohomology of a Shimura variety $X_U$ associated to a unitary
group $G$ and level subgroup $U$, with $U_p$ a {\em parahoric} subgroup of 
$G(\QQ_p)$.   Given such a $\overline{\rho}$, an analogue of Mazur's principle
should give conditions under which $\overline{\rho}$ arises
from a Shimura variety with ``less level structure at $p$''; i.e.,
under which $\overline{\rho}$ arises in the cohomology of some Shimura
variety $X_{U^{\prime}}$, with $U^{\prime}$ a level subgroup that
is equal to $U$ at primes away from $p$, but such that $U^{\prime}_p$
properly contains $U_p$.

Let us be more precise.  Let $G$ be a unitary group over $\QQ$,
such that $G(\infty) \cong \GU(1,n-1)$,
and let $p$ be a prime such that $G(\QQ_p) \cong \GL_n(\QQ_p)$.  (We will
consider a somewhat more general situation in the body of the paper; see
section~\ref{sec:basic} for the exact hypotheses.)  Up to
conjugacy, the parahoric subgroups are those subgroups of
$\GL_n(\ZZ_p)$ consisting of matrices which are ``block upper triangular''
modulo $p$.  In particular if $i_0, \dots, i_r$ are positive integers with
$i_0 + \dots + i_r = n$, let $\Gamma_{i_0, \dots, i_r}$ denote
the subgroup of $\GL_n(\ZZ_p)$ consisting of matrices which have
block form:
$$
\begin{pmatrix}
A_0 & * & * & \dots & *\\
0 & A_1 & * & \dots & *\\
& & \vdots & & \\
0 & 0 & \dots & 0 & A_r\\
\end{pmatrix}
$$
modulo $p$, where $A_j$ an $i_j$ by $i_j$ matrix.
Then the parahoric subgroups of $\GL_n(\QQ_p)$ are those subgroups conjugate 
in $\GL_n(\QQ_p)$ to some $\Gamma_{i_0, \dots, i_r}$.  The smallest of these
is the Iwahori subgroup $\Gamma_{1, \dots, 1}$; the largest is
$\Gamma_n = \GL_n(\ZZ_p)$.  Note that (even up to conjugacy) the set
of parahoric subgroups is only partially ordered by inclusion.

Now, given a mod $l$ representation $\overline{\rho}$ arising from a 
Shimura variety
$X_U$, coming from a compact open subgroup $U$ of $G(\AA^{\infty}_{\QQ})$ with
$U_p$ conjugate to $\Gamma_{i_0, \dots, i_r}$, we can ask for both necessary 
and sufficient conditions under which $\overline{\rho}$ arises from level 
$U^{\prime}$,
where $U^{\prime}$ is equal to $U$ away from $p$ but $U^{\prime}_p$
is a different parahoric subgroup conjugate to some
$\Gamma_{i_0^{\prime}, \dots, i^{\prime}_{r^{\prime}}}$ containing
$\Gamma_{i_0, \dots, i_r}$.

In section~\ref{sec:weight} we give a necessary criterion.  In particular,
we show that if $\overline{\rho}$ arises from a level whose $p$-part is 
$\Gamma_{i_0, \dots, i_r}$, then the monodromy operator $N$ on the
restriction of $\overline{\rho}$ to a decomposition group at $p$ satisfies
$N^{r+1} = 0$.  In particular $\overline{\rho}$ cannot arise from
$\Gamma_{i_0^{\prime}, \dots, i^{\prime}_{r^{\prime}}}$ with $r^{\prime} < r$
unless $N^r$ is also zero.  We then conjecture (under suitable
hypotheses on $\overline{\rho}$; c.f. Conjecture~\ref{conj:nilp}) that the
converse also holds, i.e. that if $\overline{\rho}$ satisfies certain technical
conditions, and $N^{r^{\prime}} = 0$ on $\overline{\rho}$ for some 
$r^{\prime} < r$,
then $\overline{\rho}$ arises from a level $U^{\prime}$ isomorphic to $U$ away
from $p$ but conjugate to some 
$\Gamma_{i_0^{\prime}, \dots, i^{\prime}_{r^{\prime}}}$ at $p$.

Our main results (Theorems~\ref{thm:main1} and~\ref{thm:main2}) establish 
many cases of this conjecture for $n=3$.  The techniques used are, in
spirit at least, close to those used to establish Mazur's principle.  In
particular, the two key tools used in the proof of Mazur's principle are
the Deligne-Rapoport model for elliptic curves with $\Gamma_0(p)$-level
structure and the reduction theory of Jacobians.  We study models of
the Shimura varieties under consideration in section~\ref{sec:reduction}. 
As Jacobians are unavailable for Shimura varieties of dimension 
greater than one, we instead use the Rapoport-Zink weight spectral 
sequence~\cite{RZ2} to compare the cohomology of the Shimura varieties
under consideration in characteristics $0$ and $p$.  This provides a
higher-dimensional analogue of the reduction theory of Jacobians of curves. 

These results suggest several interesting questions that are beyond the
scope of this paper.  One can ask, for instance, what happens for $n>3$.
In this case the situation seems much more complicated.  In particular,
the key point in the argument is to show that the weight
spectral sequence degenerates at $E_1$ when localized at the maximal
ideal of the Hecke algebra that corresponds to the Galois representation
being considered.  To prove this when $n=3$ we exploit a peculiarity
of the geometry in this case.  In particular, let $U$ be a subgroup of
$G(\AA^{\infty}_{\QQ})$ that is conjugate to $\Gamma_{1,2}$ or $\Gamma_{1,1,1}$
at $p$, and let $U^{\prime}$ be equal to $U$ away from $p$ and maximal
compact at $p$.  Then if $X_U$ and $X_{U^{\prime}}$ are the
corresponding Shimura varieties, we show (Propositions~\ref{prop:max1}
and~\ref{prop:iwahori}) that the intersection of any two distinct
irreducible components of the special fiber at $p$ of $X_U$ is
isomorphic to a subvariety of $X_{U^{\prime}}$.  This is crucial
in our analysis of the weight spectral sequence attached to $X_U$.
When $n$ is larger than $3$, this phenomenon only occurs for
a few level structures at $p$; essentially only those of
the form $\Gamma_{1,n-1}$ or $\Gamma_{n-1,1}$ at $p$.  The techniques
described here thus apply to these cases as well, but handling any other
level structure for $n>3$ appears to require some additional insight.

One can also ask, given a representation $\overline{\rho}$, exactly which level
structures at $p$ give rise to it.  Conjecture~\ref{conj:nilp} gives
a complete answer to this when $n=3$ (since the subgroups $\Gamma_{1,2}$
and $\Gamma_{2,1}$ are conjugate and this give rise to the same set of
representations), but not for any larger $n$.  Giving even a conjectural
answer to this question would involve a closer study of the reductions
of Shimura varieties with these level groups, and the implications of
this reduction theory for the associated Galois representations.

Finally, one can ask the same question for $p$ such that
$G(\QQ_p)$ is not isomorphic to $\GL_n$.  Again, the main difficulty seems
to be that little is known about the reduction of unitary Shimura varieties
with parahoric level structure at such $p$. 

%%%%%%%%%%%%%%%%%%%%%%%%%%%%%%%%%%%%%%%%%%%%%%%%%%%%%%%%%%%%%%%%%%%%%%%%%
\section{Unitary Shimura Varieties} \label{sec:basic}

We begin with the definition and basic properties of unitary
Shimura varieties.

Fix a totally real field $F^+$, of degree $d$ over $\QQ$.  Let $E$ be an
imaginary quadratic extension of $\QQ$, and let $x$ be
a purely imaginary element of $E$.  Let $F$ be the field 
$EF^+$.

Fix a complex square root $x_{\CC}$ of the rational number $x^2$.  
Then any embedding $\tau: F^+ \rightarrow \RR$ induces two embeddings 
$p_{\tau}, q_{\tau}: F \rightarrow \CC$, via
\begin{eqnarray*}
p_{\tau}(a+bx) & = & \tau(a) + \tau(b)x_{\CC} \\
q_{\tau}(a+bx) & = & \tau(a) - \tau(b)x_{\CC}.
\end{eqnarray*}

Fix an integer $n$, and a central simple algebra $D$ over $F$ of $F$-dimension
$n^2$, such that $D$ is either split or a division algebra at every place
of $F$.
Also fix an involution of the second kind $\alpha \mapsto \alpha^*$ on $D$.  
Let $\CV$ be a free left $D$-module of rank one, equipped
with an alternating, nondegenerate
pairing $\langle, \rangle: \CV \times \CV \rightarrow \QQ.$  We require that
$$\langle \alpha x,y\rangle = \langle x, \alpha^* y \rangle$$ 
for all $\alpha$ in $F$.

Note that $D$ is split at infinity, as $F$ is a CM-field.
Thus for each real embedding $\tau$ of $F^+$, $D \otimes_{F^+,\tau} \RR$
is isomorphic to $M_n(\CC)$.  For a suitable choice of idempotent $e$
in $D \otimes_{F^+,\tau} \RR$, satisfying $e = e^*$, 
$e(\CV \otimes_{F^+,\tau} \RR)$ is an $n$-dimensional complex vector space,
on which $\langle, \rangle$ induces a nondegenerate, Hermitian pairing.
We let $r_{\tau}(\CV)$ and $s_{\tau}(\CV)$ denote the number of $1$'s and
$-1$'s in the signature of this pairing, respectively; this is independent
of $e$ and satisfies $r_{\tau} + s_{\tau} = n$.

We choose a particular real embedding $\tau_0$ of $F^+$, and demand
that $r_{\tau_0} = 1$ (so $s_{\tau_0} = n-1$), and that
$r_{\tau} = 0$ for all $\tau \neq \tau_0$ (so that $s_{\tau} = n$
for such $\tau$).  (This is effectively a condition on the involution
$\alpha \mapsto \alpha^*$ chosen above.)

Let $G$ be the algebraic group over $\QQ$ such that for any $\QQ$-algebra
$R$, $G(R)$ is the subgroup of $\Aut_D(\CV \otimes_{\QQ} R)$ consisting
of all $g$ such that there exists an $r$ in $R^{\times}$ with
$\langle gx,gy \rangle = r\langle x,y \rangle$ for all $x$ and $y$ in
$V \otimes_{\QQ} R$.  The discussion in the previous paragraph shows
that $G(\RR)$ is the subgroup of  
$$\prod_{\tau: F^+ \rightarrow \RR} \GU(r_{\tau}, s_{\tau})$$ 
consisting of those $g$ for which the ``similitude ratio'' of
$g_{\tau}$ is the same for each $\tau$.

Now fix a compact open subgroup $U$ of $G(\AA^{\infty}_{\QQ})$,
and consider the Shimura variety associated to $G$ and $U$.  This
is a smooth variety, defined over the field $F$, which we denote by
$(X_U)_F$. Its dimension is given by the formula 
$$\dim (X_U)_F = \sum_{\tau: F^+ \rightarrow \RR} r_{\tau}s_{\tau} = n-1.$$  
If $U$ is sufficiently small, this variety can be thought of as a fine 
moduli space for abelian varieties with PEL structures.  
{\em We assume henceforth that this is the case.}  (Note that
this is not merely a convenience to avoid the use of stacks; even in the
case of classical modular forms, Mazur's principle does not quite work
as advertised if one allows arbitrary level structures and considers
modular curves that are only coarse moduli spaces; c.f. 
Remark~\ref{remark:compgrp}.)

In order to describe this moduli problem explicitly, we first note that
if $A$ is an abelian variety (up to isogeny) over an $F$-scheme $S$, with 
an action
of $D$, then the $F$-vector space $\Lie(A/S)$ decomposes as a direct
sum
$$\Lie(A/S) = \Lie(A/S)^+ \oplus \Lie(A/S)^-,$$
where $E \subset D$ acts on $\Lie(A/S)^+$ via the inclusion
$E \subset F$ and on $\Lie(A/S)^-$ by its complex conjugate.

The scheme $(X_U)_F$ then represents the functor that associates to
an $F$-scheme $S$ the set of isomorphism classes of triples 
$(A,\lambda,\psi)$ where:
\begin{enumerate}
\item $A$ is an abelian scheme (up to isogeny) over $S$ of dimension $n^2d$, 
with an action of $D$,
\item The space $\Lie(A/S)^+$ is a locally free $\OO_S$-module of rank
$n$, and $F^+ \subset D$ acts on $\Lie(A/S)^+$ via the natural
inclusion of $F^+$ in $F$.  (Harris and Taylor~\cite{HT} refer to such
an action of $D$ on $A$ as {\em compatible}.)
\item $\lambda$ is a polarization of $A$, such
that the Rosati involution associated to $\lambda$ induces the involution 
$\alpha \mapsto \alpha^*$ on $D \subset \End_0(A)$.
\item $\psi$ is a $U$-orbit of isomorphisms 
$\CV \otimes \AA^{\infty}_{\QQ} \rightarrow V A$, sending the Weil pairing on 
$\CV \otimes \AA^{\infty}_{\QQ}$
to a scalar multiple of the pairing $\langle, \rangle$ on $VA$.  (Here
$V A$ denotes the adelic Tate module of $A$.)
\end{enumerate}

Fix a prime $p$ split in $E$, and unramified in $F^+$,
such that $D$ is split at all primes of $F$ above $p$. 
Choose a prime $\pp$ of $F$ over $p$, and let 
$\overline{\pp} = \pp \cap \OO_{F^+}$.  We wish to study models
for $(X_U)_F$ over $\OO_{F,\pp}$.  This entails some subtleties involving the
$p$-part of the level structure.  We more or less follow~\cite{TY}, section 2.

First, we need to make some assumptions on
$U$.  We first assume that $U = U_pU^{(p)}$, where $U_p$ is a compact
open subgroup of $G(\QQ_p)$ and $U^{(p)}$ is a compact open subgroup
of $G(\AA^{\infty,p}_{\QQ})$.

Observe that, since $p$ splits in $E$, we have
$$G(\QQ_p) \cong \QQ_p^{\times} \times \prod_{\overline{\pp}^{\prime} | p}
\GL_n(F^+_{\overline{\pp}^{\prime}}),$$  
where the product is over the primes $\overline{\pp}^{\prime}$ of
$F^+$ over $p$.
If we fix a maximal order $\OO$ of $D$, stable under $\alpha \mapsto \alpha^*$,
we can choose this isomorphism in such a way that $\OO \otimes \ZZ_p$
is identified with the subgroup
$$\ZZ_p^{\times} \times \prod_{\overline{\pp}^{\prime} | p}
\GL_n(\OO_{F^+,\overline{\pp}^{\prime}}).$$
Fix such an isomorphism, and
assume that $U_p$ is given by the product
$$\ZZ_p^{\times} \times \prod_{\overline{\pp}^{\prime} | p} 
U_{\overline{\pp}^{\prime}},$$
where each $U_{\overline{\pp}^{\prime}}$ is a compact open subgroup
of $\GL_n(\OO_{F^+,\overline{\pp}^{\prime}})$.

Let $i_0, \dots, i_r$ be a sequence of positive integers satisfying
$$ \sum_j i_j = n.$$  We let $\Gamma_{i_0, \dots, i_r}$ denote the
subgroup of $\GL_n(\OO_{F^+,\overline{\pp}})$ consisting of those
matrices which are block upper triangular of the form:
$$
\begin{pmatrix}
A_1 & * & * & \dots & *\\
0 & A_2 & * & \dots & *\\
& & \vdots & & \\
0 & 0 & \dots & 0 & A_r\\
\end{pmatrix}
$$
modulo $\overline{\pp}$, where $A_j$ is an $i_j$ by $i_j$ matrix for each
$j$.  (Thus if $r = 0$, we have $i_0 = n$ and $\Gamma_{i_0, \dots, i_r}$
is all of $\GL_n(\OO_{F^+,\overline{\pp}})$.)

Now assume $U_{\overline{\pp}}$ is equal to $\Gamma_{i_0, \dots, i_r}$
for some $i_0, \dots, i_r$, and let $(A,\lambda)$ be a pair over $S$ as 
above.  If $S$ is an $F$-scheme, then giving a $U_{\overline{\pp}}$-orbit
of isomorphisms $T_{\overline{\pp}} \rightarrow (T_p A)_{\overline{\pp}}$
is equivalent to giving an $\OO$-stable tower of subgroup schemes
$$0 = G_0 \subset G_1 \subset \dots \subset G_{r+1} = A[\pp],$$
where $G_{j+1}/G_j$ is a finite flat group scheme over $S$ of order
$(\# \OO/\pp \OO)^{i_j}$.

We thus define a $U_{\overline{\pp}}$-level structure on a suitable 
$(A,\lambda)$ over an $\OO_{F,\pp}$-scheme $S$ to be a tower of subgroups
$$0 = G_0 \subset G_1 \subset \dots \subset G_{r+1} = A[\pp]$$
satisfying the conditions
\begin{enumerate}
\item $G_{j+1}/G_j$ is a finite flat group scheme of order
$(\# \OO/\pp \OO)^{i_j}$ for all $j$, and
\item The cokernel of the map 
$$H^1_{\DR}((A/G_{j+1})/S) \rightarrow H^1_{\DR}((A/G_j)/S)$$
induced by the natural map $A/G_j \rightarrow A/G_{j+1}$ is a locally
free $\OO_S/\pp \otimes_{\OO_F} \OO$-module of rank $i_j$
for all $j$.
\end{enumerate}
(The second of these conditions is vacuous in characteristic zero, but
is necessary in mixed characteristic to prevent the appearance of components
contained entirely in the special fiber of the moduli space.)

Following~\cite{TY}, who only discuss the case where
$U_{\overline{\pp}} = \Gamma_{1,1,\dots,1}$, we define 
$(X_U)_{\OO_{F,\pp}}$ to be the moduli space parametrizing
tuples $(A,\lambda,\psi^{(p)},\{\psi_{\pp^{\prime}}\}, G_1, \dots, G_r)$, 
where:

\begin{enumerate}
\item $A$ is an abelian scheme over $S$ up to prime-to-$p$ isogeny, with
an action of $\OO$, such that $\OO_F \subset \OO$ acts on $\Lie(A/S)^+$
via the natural inclusion of $\OO_F$ in $\OO_{F,\pp}$,
\item $\lambda$ is a prime-to-$p$ polarization, such that the associated
Rosati involution induces $\alpha \mapsto \alpha^*$ on $\OO$, 
\item $\psi^{(p)}$ is a $U$-orbit of isomorphisms $\CV \otimes \AA^{f,p}_{\QQ}
\rightarrow V^p A$, where $V^p A$ is the prime-to-$p$ adelic Tate module
of $A$,
\item for each prime $\pp^{\prime}$ of $F$ over $\pp$ such that
$\pp^{\prime} \neq \pp$, but $\pp^{\prime} \cap \OO_E = \pp \cap \OO_E$, 
$\psi_{\pp^{\prime}}$ is a $U_{\overline{\pp}^{\prime}}$-orbit of isomorphisms
$(T_p A)_{\pp^{\prime}} \cong \OO_{\pp^{\prime}}$ (Note that it is exactly
for such $\pp^{\prime}$ that such an isomorphism exists.)
\item $G_1, \dots, G_r$ is a $U_{\overline{\pp}}$-level structure in
the above sense.  (Our convention is that $G_0$ and
$G_{r+1}$ are fixed and equal to $0$ and $A[\pp]$ respectively, and thus
do not need to appear in the moduli data.)
\end{enumerate}

The generic fiber of this scheme is isomorphic to the base change of 
$(X_U)_F$ to $F_{\pp}$.

We conclude this section by observing that the matrix
$$g_{i_1} = 
\begin{pmatrix} 0 & I_{i_2 + \dots + i_r} \\ pI_{i_1} & 0 \end{pmatrix}$$
satisfies 
$g_{i_1}\Gamma_{i_1, \dots, i_r}g_{i_1}^{-1} = \Gamma_{i_2, \dots, i_r, i_1}.$
If $U$ and $U^{\prime}$ are compact open subgroups of $G(\AA^{\infty}_{\QQ})$
of the form considered above, isomorphic except at $\pp$ and
satisfying $U_{\pp} = \Gamma_{i_1, \dots, i_r}$ and
$U^{\prime}_{\pp} = \Gamma_{i_2, \dots, i_r, i_1}$, then
$g_{i_1}$ induces an isomorphism of $(X_U)_{\OO_F/\pp}$ with
$(X_{U^{\prime}})_{\OO_F/\pp}$.  This isomorphism takes a tuple
$(A,\lambda,\psi^{(p)},\{\psi_{\pp^{\prime}}\}, G_1, \dots, G_r)$
to the tuple 
$(A/(G_1 + G_1^{\perp}),\lambda^{\prime},f \circ \psi^{(p)},
\{f \circ \psi_{\pp^{\prime}}\}, G_2/G_1, \dots, G_r/G_1, A[\pp]/G_1)$,
where $G_1^{\perp}$ denotes the subgroup of $A[\pp^*]$ that
annihilates $G_1$ under the Weil pairing
$A[\pp] \times A[\pp^*] \rightarrow \mu_p$.

%%%%%%%%%%%%%%%%%%%%%%%%%%%%%%%%%%%%%%%%%%%%%%%%%%%%%%%%%%%%%%%%%%%%%%%%%
\section{Reduction of unitary Shimura varieties} 
\label{sec:reduction}

We now discuss the mod $\pp$ reduction of $(X_U)_{\OO_{F,\pp}}$,
for $U$ such that $U_{\overline{\pp}} = \Gamma_{i_0, \dots, i_r}$
for some $i_0, \dots, i_r$.  We begin by discussing a local model
for $(X_U)_{\OO_{F,\pp}}$, first constructed by Rapoport-Zink~\cite{RZ1}.

For $j$ between $0$ and $r$, let $M_j$ be a free $\OO_{F,\pp}$-module
of rank $n$, with basis $e^j_1, \dots, e^j_n$. Define 
$$f_j: M_{j+1} \rightarrow M_j$$
by $f_j(e^{j+1}_k) = p e^j_k$ if 
$$\sum_{m=0}^{j-1} i_m < k \leq \sum_{m=0}^j i_m$$ 
and $f_j(e^j_k) = e^{j+1}_k$ otherwise.
Also define
$$f_r: M_0 \rightarrow M_r$$
by $f_r(e^0_k) = p e^r_k$ if $k > \sum_{m=0}^{r-1} i_m$
and $f_r(e^0_k) = e^r_k$ otherwise. 
Note that $f_r \circ f_{r-1} \circ \dots \circ f_0$ is equal to
multiplication by $p$.

Let $\CM$ be the $\OO_{F,\pp}$-scheme parametrizing families
$(V_0, \dots, V_r)$, where: 
\begin{itemize}
\item each $V_j$ is a local direct summand of $M_j$ that is locally free of
rank one,
\item $f_j(V_j) \subset V_{j+1}$ for $0 \leq j \leq r-1$, and
\item $f_r(V_r) \subset V_0$.  
\end{itemize}

Then $\CM$ is a regular, flat
$\OO_{F,\pp}$-scheme of relative dimension $n-1$.  Its special fiber is
a divisor with normal crossings.
More precisely, let $\CM_j$ denote the subscheme of $\CM$
parametrizing those $(V_0,\dots,V_r)$ for which $f_j(V_{j+1}) = 0$.
(Of course, we take $\CM_r$ to be the subscheme of $\CM$ for which
$f_r(V_0) = 0$.)
Then a straightforward calculation shows that $\CM_j$ is a smooth subvariety 
of the special fiber of $\CM$, and that this special fiber is simply the 
union of the subvarieties $\CM_j$. 

Following Rapoport-Zink, we will show that $\CM$ is a local model for 
$(X_U)_{\OO_{F,\pp}}$.  Let the universal object for
$(X_U)_{\OO_{F,\pp}}$ be $(\CA, \lambda^{\univ}, \psi^{(p)}, 
\{\psi_{\overline{\pp}^{\prime}}\}, \CG_1, \dots, \CG_{r})$.
Set $\CA_j = \CA/\CG_j$, and let $\phi_j$ denote the
natural map $\CA_j \rightarrow \CA_{j+1}$ for each $j$.

Fix an idempotent $e$ in $\OO_{\pp}$ such that $e = e^*$
and such that $e(\OO_{\pp})$ is a free $\OO_{F,\pp}$-module of rank $n$.
Also define, for any $\OO_{F,\pp}$-scheme $S$, and any 
$\OO \otimes_{\ZZ} \OO_F$-module $N$, the submodule $N_0$
to be the largest submodule of $N$ on which $\OO_F \subset \OO$
acts via the natural map $\iota: \OO_F \rightarrow \OO_S$.

Then we can define the scheme $\tX_U$ to be the $(X_U)_{\OO_{F,\pp}}$-scheme
that represents the functor taking an $(X_U)_{\OO_{F,\pp}}$ scheme
$\theta: S \rightarrow (X_U)_{\OO_{F,\pp}}$ to the set of diagrams:
$$
\begin{array}{ccccccc}
eH^1_{\DR}(\theta^*\CA_0/S)_0 & \stackrel{\phi_0}{\leftarrow} &
eH^1_{\DR}(\theta^*\CA_1/S)_0 & \stackrel{\phi_1}{\leftarrow} & \dots &
\stackrel{\phi_r}{\leftarrow} & eH^1_{\DR}(\theta^*\CA_{r+1}/S)_0 \cong
eH^1_{\DR}(\theta^*\CA_0/S)_0 \\
\downarrow & & \downarrow & & & & \downarrow\\
(M_0)_S & \stackrel{f_0}{\leftarrow} & (M_1)_S & \stackrel{f_1}{\leftarrow} &
\dots & \stackrel{f_r}{\leftarrow} & (M_0)_S\\
\end{array}
$$
where the vertical arrows are isomorphisms and the leftmost and rightmost
vertical arrow ``coincide''.

We have an obvious map of $\tX_U$ to $(X_U)_{\OO_{F,\pp}}$.  We also
have a map of $\tX_U$ to $\CM$, which associates to any 
$\theta: S \rightarrow (X_U)_{\OO_{F,\pp}}$ and any diagram as above
the tuple $((V_0)_S, \dots, (V_r)_S)$, where $V_j$
corresponds to the subbundle $e\Lie(\theta^*\CA_j/S)_{\pp}^*$ of
$eH^1_{\DR}(\theta^*\CA_j/S)_{\pp}$ under the given identification
of $eH^1_{\DR}(\theta^*\CA_j/S)_{\pp}$ with $(M_j)_S$.
The methods of Rapoport-Zink then show that
$$(X_U)_{\OO_{F,\pp}} \leftarrow \tX_U \rightarrow \CM$$
is a local model diagram; i.e., that both maps are smooth of the same
dimension.

In particular $\CM$ is a local model for $(X_U)_{\OO_{F,\pp}}$, so 
$(X_U)_{\OO_{F,\pp}}$ is regular and flat over $\OO_{F,\pp}$,
and its special fiber $(X_U)_{\OO_F/\pp}$ is a divisor with normal crossings.

\begin{definition} Let $(X_U)_j$ denote the subscheme of
$(X_U)_{\OO_F/\pp}$ on which 
$$\phi_j: \Lie(\CA_j/S) \rightarrow \Lie(\CA_{j+1}/S)$$
is the zero map.
\end{definition}

The base change of $(X_U)_j$ to $\tX_U$
coincides with the base change of $\CM_j$ to $\tX_U$.  It follows that
$\CM_j$ is a local model for $(X_U)_j$ for all $j$; i.e., that the
$(X_U)_j$ are smooth $\OO_F/\pp$-schemes of dimension $n-1$.  Their union 
is the special fiber of $(X_U)_{\OO_{F,\pp}}$. 

Our goal is to study the global geometry of the schemes
$(X_U)_j$ for various $U$.  In general this is quite complicated,
so we will consider only a few special cases here.  Fix a $U$ such that 
$U_{\overline{\pp}} = \Gamma_n$; i.e., is maximal compact.  Let
$U_{1,n-1}$ be the subgroup of $U$ such that 
$(U_{1,n-1})_{\overline{\pp}} = \Gamma_{1,n-1}$ and 
$U$ is isomorphic to $U_{1,n-1}$
``away from $\overline{\pp}$''.  Then 
$(X_U)_{\OO_F/\pp}$ is smooth, whereas 
$(X_{U_{1,n-1}})_{\OO_F/\pp}$ breaks up into two pieces,
$(X_{U_{1,n-1}})_0$ and $(X_{U_{1,n-1}})_1$.  

Understanding these pieces will require an analysis of the $\pp$-torsion
of the moduli objects for $(X_U)_{\OO_F/\pp}$.  We approach this by
Dieudonn{\'e} theory.  In particular, let $k$ be an
algebraically closed field of characteristic $p$, and
let $(A,\lambda,\psi^{(p)},\{\psi_{\pp^{\prime}}\})$
be a $k$-point of $(X_U)_{\OO_F/\pp}$.  Let $D_{\pp} A$
denote the contravariant Dieudonn{\'e} module of $A[\pp^{\infty}]$.
It is a free $W(k)$-module of rank $n^2d_{\pp}$, where $d_{\pp}$
is the residue class degree of $\pp$ over $p$.  It is equipped with
the following extra structures:
\begin{itemize}
\item a $\sigma$-linear endomorphism $F: D_{\pp} A \rightarrow D_{\pp} A$,
where $\sigma$ is the Witt vector Frobenius on $W(k)$,
\item a $\sigma$-antilinear endomorphism $V: D_{\pp} A \rightarrow D_{\pp} A$,
such that $FV = VF = p$, and
\item an action of $\OO \otimes_{\ZZ} W(k)$, induced by the action of $\OO$
on $A$, that commutes with $F$ and $V$.
\end{itemize}
 
The map $\spec k \rightarrow (X_U)_{\OO_{F,\pp}}$ corresponding to
the point $(A,\lambda,\psi^{(p)},\{\psi_{\pp^{\prime}}\})$
yields a map $\spec k \rightarrow \spec \OO_{F,\pp}$, and hence a map
$\OO_F \rightarrow k$.  This in turn induces a map
$\iota: \OO_F \rightarrow W(\overline{\FF}_p)$.  For $i$ in
$\ZZ/d_{\pp}\ZZ$ 
let $(D_{\pp} A)_i$ denote the submodule of $D_{\pp} A$
on which $\OO_F \subset \OO$ acts via $\sigma^i \iota: \OO_F \rightarrow
W(\overline{\FF}_p)$.  (Here $\sigma$ is the Witt vector Frobenius.)
Then $D_{\pp} A$ breaks up as the direct sum of $(D_{\pp} A)_i$ for
$i \in \ZZ/d_{\pp}\ZZ$.  Each $(D_{\pp} A)_i$ has rank $n^2$ over
$W(k)$, and comes equipped with an action of $\OO \otimes_{\OO_F} W(k)$.
The latter is isomorphic to $\GL_n(W(k))$.

Note that the semilinearity properties of $F$ and $V$ imply that they induce
maps
$$F: (D_{\pp} A)_i \rightarrow (D_{\pp} A)_{i+1}$$
$$V: (D_{\pp} A)_{i+1} \rightarrow (D_{\pp} A)_i$$
for all $i$.

Let $\overline{D}_{\pp} A$ denote $D_{\pp} A/p D_{\pp} A$.  Then
(see for instance~\cite{oda}), we have an isomorphism
$H^1_{\DR}(A/k)_{\pp} \cong \overline{D}_{\pp} A$.  Moreover, the submodule
$V\overline{D}_{\pp} A$ corresponds under this isomorphism to
the submodule $\Lie(A/k)^*_{\pp}$ of $H^1_{\DR}(A/k)_{\pp}$.  Since
(by our definition of the moduli problem)
$\OO_F$ acts on $\Lie(A/k)^*_{\pp}$ via the natural map $\OO_F \rightarrow k$,
it follows that $V\overline{D}_{\pp} A$ lies in $(\overline{D}_{\pp} A)_0$.
Moreover, it has dimension $n$ as a $k$-vector space.  Thus the maps
$$V: (\overline{D}_{\pp} A)_{i+1} \rightarrow (\overline{D}_{\pp} A)_i$$
are zero for $i \neq 0$, and the map
$$V: (\overline{D}_{\pp} A)_1 \rightarrow (\overline{D}_{\pp} A)_0$$
has an $n$-dimensional image.

On the other hand, the kernel of $F$ on $\overline{D}_{\pp} A$ is
equal to the image of $V$.  Thus the maps
$$F: (\overline{D}_{\pp} A)_i \rightarrow (\overline{D}_{\pp} A)_{i+1}$$
are isomorphisms for $i \neq 0$, and the map
$$V: (\overline{D}_{\pp} A)_0 \rightarrow (\overline{D}_{\pp} A)_1$$
has an $n^2-n$-dimensional image.

Consider the endomorphism $F^{d_{\pp}}$ of
the Dieudonn{\'e} module $\overline{D}_{\pp} A$.  It induces a map
$$(\overline{D}_{\pp} A)_0 \rightarrow (\overline{D}_{\pp} A)_0$$
whose image has dimension $n^2-n$.  The kernel of this map is
$V((\overline{D}_{\pp} A)_1)$.  Since both of these are stable under
$\OO$, either the kernel is contained in the image or intersects it
trivially.  In the former case $F^{id_{\pp}}(\overline{D}_{\pp} A)_0$
has dimension $n^2 - n$ for all $i \geq 1$; in the latter case
$F^{2d_{\pp}}(\overline{D}_{\pp} A)_0$ has dimension $n^2 - 2n$.
Once again, this space will either contain the kernel of
$F^{d_{\pp}}$ or meet it trivially.

Proceeding in this fashion, we see that there is a unique $h$ between $0$
and $n - 1$, such that $F^{id_{\pp}}(\overline{D}_{\pp} A)_0$ has
dimension $nh$ for all $i \geq n - h$.  It follows that, in the language
of~\cite{HT}, Chapter II, $eA[\pp^{\infty}]$ is a Barsotti-Tate 
$\OO_F$-module of {\'e}tale height $h$.  Following~\cite{HT}, we
set $(X_U)^{(h)}$ to be the locus of points in $(X_U)_{\OO_F/\pp}$
with this property, and let $(X_U)^{[h]}$ denote its closure.
Harris and Taylor show (Corollary III.4.4 of~\cite{HT}) that 
$(X_U)^{(h)}$ is smooth of pure dimension $h$.

Let $(\overline{D}_{\pp} A)^{et}$ denote the submodule
$F^{n-h}(\overline{D}_{\pp} A)$ of $\overline{D}_{\pp} A$.  It is the
largest Dieudonn{\'e} submodule of $\overline{D}_{\pp} A$ on which $F$ is
an isomorphism.  Let $(\overline{D}_{\pp} A)^{\conn}$ be the submodule
$$\{x \in \overline{D}_{\pp} A: \exists i: F^ix = 0\}$$
of $(\overline{D}_{\pp} A)$.  The above analysis shows that
we have a decomposition
$$\overline{D}_{\pp} A = (\overline{D}_{\pp} A)^{\conn} \oplus
(\overline{D}_{\pp} A)^{\et}.$$

With these facts established, we may now show:

\begin{proposition} \label{prop:max0} The natural map 
$$(X_{U_{1,n-1}})_{\OO_F/\pp} \rightarrow (X_U)_{\OO_F/\pp}$$
induces an isomorphism of $(X_{U_{1,n-1}})_0$ with 
$(X_U)_{\OO_F/\pp}$.
\end{proposition}
\begin{proof}
This natural map is clearly proper, so it suffices to show that it is
a bijection on $\overline{\FF}_p$-points and tangent spaces.  Let 
$x = (A,\lambda,\psi^{(p)},\{\psi_{\pp^{\prime}}\},G_1)$
be an $\overline{\FF}_p$-point of $(X_{U_{1,n-1}})_0$.  Take 
$A_1 = A/G_1$, and let $\phi_0$ be the natural map $A \rightarrow A_1$. 

Since $x$ lies in $(X_{U_{1,n-1}})_0$, $\phi_0$ induces the
zero map $\Lie(A_1)^* \rightarrow \Lie(A)^*$.  In the language of
Dieudonn{\'e} theory, this means that $\phi$ maps
$VD_{\pp} A_1$ into $pD_{\pp} A$.  Since $V$ is injective on the
Dieudonn{\'e} module of a $p$-divisible group, this means that
$\phi$ maps $D_{\pp} A_1$ into $F D_{\pp} A$.

Condition (2) in the definition of $U_{\overline{\pp}}$-level structure,
and the isomorphism between $\overline{D}_{\pp} A$
and $H^1_{\DR}(A)_{\pp}$, together imply that 
the subspace $\phi_0(\overline{D}_{\pp} A_1)_i$ of $(\overline{D}_{\pp} A)_i$ 
has dimension $n^2-n$.  Since this image is contained in 
$F(\overline{D}_{\pp} A)_{i-1}$, we find in particular that 
$\phi_0(\overline{D}_{\pp} A_1)_1$ 
is contained in $F(\overline{D}_{\pp} A)_0$; since both have dimension
$n^2-n$, the two must be equal.  But since $\phi_0(\overline{D}_{\pp} A_1)$ 
is a Dieudonn{\'e} submodule of $\overline{D}_{\pp} A$, it follows
that $\phi_0(\overline{D}_{\pp} A_1)_i$ contains $F^i(\overline{D}_{\pp} A)_0$
for all $i \geq 1$.  When $1 \leq i \leq d_{\pp}$, both of these
spaces have dimension $n^2-n$, so $\phi_0(\overline{D}_{\pp} A_1)_i$ equals
$F^i(\overline{D}_{\pp} A)_0$ for $1 \leq i \leq d_{\pp}$.

Thus $\phi_0(\overline{D}_{\pp} A_1)$ (and hence $G_1$)
is uniquely determined given $A$.  Conversely, for a 
point
$(A,\lambda,\psi^{(p)},\{\psi_{\pp^{\prime}}\})$ of
$(X_U)_{\OO_F,\pp}(\overline{\FF}_p)$, the submodule $M$ of 
$\overline{D}_{\pp} A$ defined by 
$$M_i = F^i (\overline{D}_{\pp} A)_0, \mbox{$1 \leq i \leq d_{\pp}$}$$
fits into an exact sequence
$$0 \rightarrow M_i \rightarrow \overline{D}_{\pp} A \rightarrow
C \rightarrow 0$$
of Dieudonn{\'e} modules.  Then $C$ is the Dieudonn{\'e} module of
some group $G_1$, and the surjection of $\overline{D}_{\pp} A$ onto
$C$ induces an injection of $G_1$ into $A$, and hence determines
a point of $(X_{U_{1,n-1}})_0$. 

Thus the map from $(X_{U_{1,n-1}})_0$ to $(X_U)_{\OO_{F,\pp}}$ is
a bijection on $\overline{\FF}_p$-points.  We now show it is a bijection
on tangent spaces.  Note that a first order deformation of the datum
$(A,\lambda,\psi^{(p)},\{\psi_{\pp^{\prime}}\},G_1)$
is equivalent to a first order deformation of the datum
$(A,\lambda,\psi{(p)},\{\psi_{\pp^{\prime}}\},A_1,\phi_0: A
\rightarrow A_1)$.

By a result of Grothendieck~\cite{montreal}, a first-order
deformation of $A$ is equivalent to a lift of the subspace
$\Lie(A)^*$ of $H^1_{\DR}(A)$ to a direct summand
of the $k[\epsilon]/\epsilon^2$-module 
$H^1_{\cris}(A)_{k[\epsilon]/\epsilon^2}$.  Such a deformation
will be compatible with $\lambda$ if this direct summand is isotropic
with respect to the natural pairing:
$$H^1_{\cris}(A)_{k[\epsilon]/\epsilon^2} \times 
H^1_{\cris}(A)_{k[\epsilon]/\epsilon^2} \rightarrow k[\epsilon]/\epsilon^2$$
induced by $\lambda$, and hence is uniquely determined by a lift of
$(\Lie(A)^*)^+$ to $H^1_{\cris}(A)^+_{k[\epsilon]/\epsilon^2}$.
This lift must also be compatible with the action of $\OO$.

To give a deformation of $(A,A_1,\phi_0: A \rightarrow A_1)$ lifting
a deformation as above, we must give a lift of $\Lie(A_1)^*$ to
$H^1_{\cris}(A_1)_{k[\epsilon]/\epsilon^2}$, compatible with the
action of $\OO$, such that the map
$\phi_0$ induces on crystalline cohomology maps this lift into the given
lift of $\Lie(A)^*$.  This deformation will be contained in
$(X_{U_{1,n-1}})_0$ if (and only if) $\phi_0$ induces the zero map
on the given lift of $(\Lie(A_1)^*)^+$. 
Since $\ker \phi_0 \subset A[\pp]$,
$\phi_0$ is an isomorphism on the minus part of the crystalline cohomology,
so we can restrict our attention to $(\Lie(A_1)^*)^+$.  A lift of this space
must have dimension $n$, and must be contained in 
$(H^1_{\cris}(A_1)_{k[\epsilon]/\epsilon^2})_0$.  On the other hand,
the cokernel (and hence also the kernel) of the map 
$H^1_{\DR}(A_1)_0 \rightarrow H^1_{\DR}(A)_0$, has dimension $n$.  Thus
the kernel of the map
$(H^1_{\cris}(A_1)_{k[\epsilon]/\epsilon^2})_0 \rightarrow
(H^1_{\cris}(A)_{k[\epsilon]/\epsilon^2})_0$ is free of rank $n$, and
since the lift of $(\Lie(A_1)^*)^+$ must be contained in this kernel, the two
are equal.

Thus, for any tangent vector to a point on
$(X_U)_{\OO_F/\pp}$, there is at most one tangent vector of
$(X_{U_{1,n-1}})_0$ lifting it.  In particular the map from
$(X_{U_{1,n-1}})_0$ to $(X_U)_{\OO_F/\pp}$ 
is injective on tangent spaces, and since the former and latter are
smooth of the same dimension, it is an isomorphism on tangent spaces as
required.
\end{proof}

We turn next to the geometry of $(X_{U_{1,n-1}})_1$, which is
somewhat more complicated.  First note that the map 
$$(X_{U_{1,n-1}})_{\OO_{F,\pp}} \rightarrow (X_U)_{\OO_{F,\pp}}$$ 
has degree $\frac{p^{nd_{\pp}} - 1} {p^{d_{\pp}} - 1}$.  
Since $(X_U)_{\OO_F/\pp}$ is a divisor of
$(X_U)_{\OO_{F,\pp}}$, and its pullback to $(X_{U_{1,n-1}})_{\OO_{F,\pp}}$
is simply $(X_{U_{1,n-1}})_{\OO_F/\pp}$, the map
$$(X_{U_{1,n-1}})_{\OO_F/\pp} \rightarrow (X_U)_{\OO_F/\pp}$$ 
has degree $\frac{p^{nd_{\pp}} - 1} {p^{d_{\pp}} - 1}$ as well.  
In light of the preceding proposition, this means that the map
$$(X_{U_{1,n-1}})_1 \rightarrow (X_U)_{\OO_F/\pp}$$
has degree
$$\frac{p^{nd_{\pp}} - 1} {p^{d_{\pp}} - 1} - 1 = 
p^{d_{\pp}}\frac{p^{(n-1)d_{\pp}} - 1} {p^{d_{\pp}} - 1}.$$

The following proposition gives the structure of this map:

\begin{proposition} \label{prop:max1}
The map $(X_{U_{1,n-1}})_1 \rightarrow (X_U)_{\OO_F/\pp}$ is
the composition of a finite, purely inseparable morphism 
$(X_{U_{1,n-1}})_1 \rightarrow (X_{U_{1,n-1}})_1^{\sep}$
of degree $p^{d_{\pp}}$ with a finite, separable morphism 
$(X_{U_{1,n-1}})_1^{\sep} \rightarrow (X_U)_{\OO_F/\pp}$
of degree
$\frac{p^{(n-1)d_{\pp}} - 1}{p^{d_{\pp}} - 1}.$
The latter morphism is {\'e}tale over $(X_U)^{(n-1)}$.  Moreover,
the preimage of $(X_U)^{[n-2]}$ under this morphism contains
$(X_{U_{1,n-1}})_0 \cap (X_{U_{1,n-1}})_1$, 
and the induced map:
$$ (X_{U_{1,n-1}})_0 \cap (X_{U_{1,n-1}})_1 \rightarrow (X_U)^{[n-2]}$$ 
is an isomorphism.

Over $(X_U)^{(n-2)}$, the map
$(X_{U_{1,n-1}})_1^{\sep} \rightarrow (X_U)_{\OO_F/\pp}$ 
is unbranched along the intersection of $(X_{U_{1,n-1}})_0
\cap (X_{U_{1,n-1}})_1$ with the preimage of $(X_U)^{(n-2)}$,
and branched to order $p^{d_{\pp}}$ elsewhere.
\end{proposition}

\begin{proof}
Let $x = (A,\lambda,\rho,\psi^{(p)},\{\psi_{\pp^{\prime}}\}, G_1)$
be a $k$-point of $(X_{U_{1,n-1}})_1$.  Let $\overline{x}$ denote
the image of $x$ in $(X_U)_{\OO_F/\pp}$.
Let $A_1$ and $A_2$ denote $A/G_1$ and $A/A[\pp]$,
and $\phi_0$ and $\phi_1$ denote the natural maps
$A \rightarrow A_1$ and $A_1 \rightarrow A_2$, respectively.
We know, by the moduli interpretation of $(X_{U_{1,n-1}})_1$, that
$\phi_1$ induces the zero map from $\Lie(A_2)^*$ to $\Lie(A_1)^*$.
In terms of Dieudonn{\'e} theory, this means that $\phi_1$ maps
$V D_{\pp} A_2$ into $p D_{\pp} A_1$.  On the other hand, $\phi_1 \circ \phi_0$
identifies $D_{\pp} A_2$ with $p D_{\pp} A$.  Thus $\phi_0(D_{\pp} A_1)$ 
contains $V D_{\pp} A$.  Conversely, if $M$ is any $\OO$-stable
Dieudonn{\'e} submodule
of $D_{\pp} A$ that contains $V D_{\pp} A$ and such that the quotient
$D_{\pp} A/M$ is a free $\OO \otimes k$-module of rank $1$, then
$M$ determines a unique subgroup $G_1$ of $A[\pp]$ such that the 
map $A/G_1 \rightarrow A/A[\pp]$ induces the zero map on the Lie algebras,
and hence a unique point of $(X_{U_{1,n-1}})_1$ over $\overline{x}$. 

Now suppose we have such an $M$, and
assume that $\overline{x}$ 
lies in $(X_U)^{(n-1)}$.  In this case 
$\overline{D}_{\pp} A$ breaks up as a direct sum
$\overline{D}_{\pp} A = (\overline{D}_{\pp} A)^{\conn} \oplus
(\overline{D}_{\pp} A)^{\et},$
where $(\overline{D}_{\pp} A)^{\conn}_i$ has dimension $n$ for all
$i$.  Let $\overline{M}$ denote the projection of $M$ to 
$\overline{D}_{\pp} A$; since $M$ contains $p\overline{D}_{\pp} A$
it is clear that $\overline{M}$ determines $M$. 

Since $\overline{M}$ contains $V\overline{D}_{\pp} A$,
$\overline{M}_0$ contains $(\overline{D}_{\pp} A)^{\conn}_0.$
Thus, under the projection of $\overline{D}_{\pp} A$ to 
$(\overline{D}_{\pp} A)^{\et}$, $\overline{M}_0$ gets sent to a submodule of
rank $n^2 - 2n$ of $(\overline{D}_{\pp} A)^{\et}_0$.  But
for $1 \leq i \leq d_{\pp} - 1$, $F^{d_{\pp} - i}$ induces
an isomorphism of $\overline{M}_i$ with $\overline{M}_0$, and 
hence $\overline{M}_i$ gets sent
to a submodule of rank $n^2 - 2n$ of $(\overline{D}_{\pp} A)^{\et}_i$.
A dimension count shows that $\overline{M}_i$ therefore contains 
$(\overline{D}_{\pp} A)^{\conn}_i$ for all $i$.  In particular
$\overline{M}$ is the preimage of a unique $\OO$-stable submodule 
$\overline{M}^{\et}$
of $(\overline{D}_{\pp} A)^{\et}$, with $\dim \overline{M}^{\et}_i = n^2 -2n$
for all $i$, and any such submodule determines a unique $M$.

Moreover, $\overline{M}^{\et}$ determines, and is uniquely determined by,
the submodule $e\overline{M}^{\et}$ of $e(\overline{D}_{\pp} A)^{\et}$.
The latter is the Dieudonn{\'e} module of the constant
group scheme $(\OO_F/\pp)^{n-1}$, and $e(\overline{D}_{\pp} A)^{\et}/eM$
is the Dieudonn{\'e} module of $(\OO_F/\pp)$.  Thus points of
$(X_{U_{1,n-1}})_1$ over $\overline{x}$ correspond to subgroups 
of $(\OO_F/\pp)^{n-1}$, isomorphic to $\OO_F/\pp$, and there are
$\frac {p^{(n-1)d_{\pp}} - 1}{p^{d_{\pp}} - 1}$ of these.  It follows
that the map $(X_{U_{1,n-1}})_1 \rightarrow (X_U)_{\OO_F/\pp}$ is
the composition of a purely inseparable map of degree $p^{d_{\pp}}$
with a separable map of degree
$\frac {p^{(n-1)d_{\pp}} - 1}{p^{d_{\pp}} - 1}$ that is {\'e}tale
over $(X_U)^{(n-1)}$.

Now suppose that $\overline{x}$ lies in $(X_U)^{(n-2)}$.  We again have
a decomposition
$$\overline{D}_{\pp} A = (\overline{D}_{\pp} A)^{\conn} \oplus
(\overline{D}_{\pp} A)^{\et},$$
but now $e(\overline{D}_{\pp} A)^{\et}$ is the Dieudonn{\'e} module
of the constant group scheme $(\OO_F/\pp)^{n-2}$.  The space
$e(\overline{D}_{\pp} A)^{\conn}$ is $2{d_{\pp}}$-dimensional, and
$F$ acts nilpotently on it with a one dimensional-kernel.  Thus in
particular $e(\overline{D}_{\pp} A)^{\conn}$ has exactly one Dieudonn{\'e}
submodule of any given dimension between zero and $2d_{\pp}$. 

Let $\overline{M}^{\et}$ be the image of $\overline{M}$ in 
$(\overline{D}_{\pp} A)^{\et}$.
The dimension of $e\overline{M}^{\et}_i$ is either $n-2$ for all $i$, 
or $n-3$ for all $i$;
in the latter case $e\overline{M}$ contains $e(\overline{D}_{\pp} A)^{\conn}$;
in the former case the intersection of $e\overline{M}$ with
$e(\overline{D}_{\pp} A)^{\conn}$ is the unique Dieudonn{\'e} submodule
of $e(\overline{D}_{\pp} A)^{\conn}$ of dimension $d_{\pp}$, and $\overline{M}$
is uniquely determined.  In the latter case $\overline{M}$ is determined by
the choice of $e\overline{M}^{\et}$, which is equivalent to the choice
of a subgroup of $(\OO_F/\pp)^{n-2}$ isomorphic to $(\OO_F/\pp)$. 

Hence over any geometric point of $(X_U)^{n-2}$, there are
$$1 + \frac{p^{(n-2)d_{\pp}} - 1}{p^{d_{\pp}} - 1}$$ geometric
points of $(X_{U_{1,n-1}})_1$.  In order to understand the branching along 
these
points, let $S$ be the spectrum of a discrete valuation ring $R$, and
let $(A,\lambda,\psi^{(p)},\{\psi_{\pp^{\prime}}\})$ be
an $S$-point of $(X_U)_{\OO_F/\pp}$, whose generic point maps into
$(X_U)^{(n-1)}$ and whose closed point maps into $(X_U)^{(n-2)}$.
Let $K$ be the algebraic closure of its field of fractions, and let $k$
be the algebraic closure of its residue field.  Then $eA_K[\pp]^{\et}$
is isomorphic to $(\OO_F/\pp)^{n-1}$, whereas $eA_k[\pp]^{\et}$ is
isomorphic to $(\OO_F/\pp)^{n-2}$. 

There are $\frac{p^{(n-1)d_{\pp}} - 1}{p^{d_{\pp}} - 1}$ $K$-points
of $(X_U)_1$ lying over $\spec K \rightarrow (X_U)_{\OO_F/\pp}$,
corresponding to the subgroups of $eA_K[\pp]^{\et}$ isomorphic to
$\OO_F/\pp$.  Under the specialization map $eA_K[\pp] \rightarrow eA_k[\pp]$,
$p^{d_{\pp}}$ such subgroups specialize to any given subgroup of
$eA_k[\pp]$ isomorphic to $\OO_F/\pp$.  The remaining subgroup specializes
to a connected subgroup of $eA_k[\pp]$.

Thus, along the preimage of $(X_U)^{(n-2)}$, the separable part of the map 
$(X_{U_{1,n-1}})_1 \rightarrow (X_U)_{\OO_F/\pp}$ is unbranched at those
points corresponding to subgroups $G_1$ that are connected, and branched
to order $p^{d_{\pp}}$ elsewhere.  On the other hand, if $G_1$ is
connected, then the map $A \rightarrow A/G_1$ is zero on some part
of $\Lie(A)^+$.  The latter has no proper $\OO$-stable subspaces, 
so this map must be zero on all of $\Lie(A)^+$.  Thus such $G_1$
determine points that lie on $(X_{U_{1,n-1}})_0$ as well as
$(X_{U_{1,n-1}})_1$.  Conversely, it is easy to check that the image
of $(X_U)^{[n-2]}$ in $(X_{U_{1,n-1}})_0$ under the inverse of the
natural isomorphism of $(X_{U_{1,n-1}})_0$ with $(X_U)_{\OO_F/\pp}$
is precisely the intersection $(X_{U_{1,n-1}})_0 \cap (X_{U_{1,n-1}})_1$.
\end{proof}

\begin{remark} \rm
Note that the isomorphism $g_1: (X_{U_{1,n-1}})_{\OO_{F,\pp}}
\cong (X_{U_{n-1,1}})_{\OO_{F,\pp}}$ identifies 
$(X_{U_{1,n-1}})_0$ with $(X_{U_{n-1,1}})_1$, and 
$(X_{U_{1,n-1}})_1$ with $(X_{U_{n-1,1}})_0$.  Thus we obtain
results about the structure of $X_{U_{n-1,1}}$ as well.
\end{remark}

We will also need to understand the case in which $n = 3$, and the
level structure at $\overline{\pp}$ comes from the Iwahori subgroup
$\Gamma_{1,1,1}$.  Let $U_{1,1,1}$ denote the subgroup of
$G(\AA^{\infty}_{\QQ})$ that is isomorphic to $\Gamma_{1,1,1}$ at
$\overline{\pp}$ and to $U$ at all other places.  Then
$(X_{U_{1,1,1}})_{\OO_{F,\pp}}$ 
parametrizes tuples $(A,\lambda,\psi^{(p)},\{\psi_{\pp^{\prime}}\},
G_1,G_2)$.  Its fiber over $\OO_F/\pp$ breaks up into three components
$(X_{U_{1,1,1}})_0$, $(X_{U_{1,1,1}})_1$ and $(X_{U_{1,1,1}})_2$.
These components parametrize the loci where the maps
$A \rightarrow A/G_1$, $A/G_1 \rightarrow A/G_2$ and $A/G_2 \rightarrow
A/A[\pp]$ (respectively) induce the zero map on the positive parts of
the appropriate Lie algebras.

The variety $(X_{U_{1,1,1}})_0$ admits a map to $(X_{U_{1,2}})_0$ by
forgetting $G_2$, and a map to $(X_{U_{2,1}})_0$ by forgetting $G_1$.
The resulting map
$$(X_{U_{1,1,1}})_0 \rightarrow (X_{U_{1,2}})_0 \times_{(X_U)_{\OO_F/\pp}}
(X_{U_{2,1}})_0$$
is clearly proper and injective on $S$-valued points for any $S$.  It
is therefore a closed immersion.  But $(X_{U_{1,2}})_0$ is isomorphic
to $(X_U)_{\OO_F/\pp}$, so the above product is simply
$(X_{U_{2,1}})_0$.  We thus have a closed immersion of $(X_{U_{1,1,1}})_0$
into $(X_{U_{2,1}})_0$. 
On the other hand, if
$(A,\lambda,\psi^{(p)},\{\psi_{\pp^{\prime}}\},G_2)$
is an $\overline{\FF}_p$-point of $(X_{U_{2,1}})_0$ lying over
$(X_U)^{(2)}$, then we have an exact sequence
$$0 \rightarrow G_2^{\conn} \rightarrow G_2 \rightarrow G_2^{\et} \rightarrow
0,$$
and it is easy to see (using the above Dieduonn{\'e} theory calculations)
that $(A,\lambda,\psi^{(p)},\{\psi_{\pp^{\prime}}\},G_2^{\conn},
G_2)$ is a point of $(X_{U_{1,1,1}})_0$ mapping to $x$.  Thus
the natural map $(X_{U_{1,1,1}})_0 \rightarrow (X_{U_{2,1}})_0$ is
an isomorphism.  It follows that $(X_{U_{1,1,1}})_0$ is
isomorphic to $(X_{U_{1,2}})_1$.

The automorphism $g_1$ of $(X_{U_{1,1,1}})_{\OO_{F,\pp}}$
cyclically permutes $(X_{U_{1,1,1}})_0$, $(X_{U_{1,1,1}})_1$, and
$(X_{U_{1,1,1}})_2$.  Thus the latter two components are
isomorphic to $(X_{U_{1,2}})_1$ as well.

Finally, we turn to the question of the intersections
$(X_{U_{1,1,1}})_i \cap (X_{U_{1,1,1}})_{i^{\prime}}$ for
$i,i^{\prime}$ in $\{0,1,2\}$.  Note that the map
$$(X_{U_{1,1,1}})_{\OO_F/\pp} \rightarrow (X_{U_{1,2}})_{\OO_F/\pp}$$
induces a map
$$(X_{U_{1,1,1}})_0 \cap (X_{U_{1,1,1}})_1 \rightarrow
(X_{U_{1,2}})_0 \cap (X_{U_{1,2}})_1.$$
that is a bijection on points.  Since the map $(X_{U_{1,1,1}})_1
\rightarrow (X_{U_{1,2}})_1$ is an isomorphism, this induced map
is as well.  Thus $(X_{U_{1,1,1}})_0 \cap (X_{U_{1,1,1}})_1$ is
naturally isomorphic to $(X_U)^{(1)}$.  The action of
$g_1$ allows us to deduce this for the other pairwise intersections as
well.

To summarize:
\begin{proposition} \label{prop:iwahori}
The three components $(X_{U_{1,1,1}})_0$, $(X_{U_{1,1,1}})_1$,
and $(X_{U_{1,1,1}})_2$ are all isomorphic to $(X_{U_{1,2}})_1$.
Their pairwise intersections are isomorphic to $(X_U)^{(1)}$.
\end{proposition}

%%%%%%%%%%%%%%%%%%%%%%%%%%%%%%%%%%%%%%%%%%%%%%%%%%%%%%%%%%%%%%%%%%%%%%%%%
\section{Galois representations} \label{sec:representations}

Henceforth, we assume that there is some place of $F$ at which $D$ is
division algebra.  It is easily seen that this implies that
$(X_U)_{\OO_F,\pp}$ is proper over $\OO_{F,\pp}$.  This avoids complications
due to endoscopy, allowing us to apply results from~\cite{HT},
and will later allow us to apply the weight spectral
sequence to the cohomology of this variety.  

Fix a compact open subgroup $U$ of $G$, let $(X_U)_{\overline{F}}$ be
the base change of $(X_U)_F$ to $\overline{F}$, and let $\TT_U$ denote the
$\ZZ_l$-subalgebra of $\End(H^{n-1}((X_U)_{\overline{F}}, \QQ_l))$
generated by the Hecke operators at primes that split in $E$, are unramified
in $F$, and do
not divide the level of $U$.  It is a finite, flat, commutative 
$\ZZ_l$-algebra.  Moreover, since the Hecke operators in $\TT_U$
commute with their adjoints with respect to Poincar{\'e} duality on
$H^{n-1}((X_U)_{\overline{F}}, \QQ_l)$, $\TT_U$ contains no nonzero
nilpotent elements.  Let $\tm$ be a minimal prime of $\TT_U$.  Then
$\tm$ determines, for each prime $q$ that splits in $E$ and does not 
divide the level of $U$, an unramified representation $\pi_{\tm,q}$ 
of $G(\QQ_q)$.

Let $X_{\overline{F}}$ denote the ind-scheme obtained by taking the
limit of $(X_U)_{\overline{F}}$ over all $U$.
The localization $H^{n-1}_{\et}((X_U)_{\overline{F}}, \QQ_l)_{\tm}$ is
given by a direct sum
$$\bigoplus_{\pi} 
H^{n-1}(X_{\overline{F}}, \QQ_l)[\pi]^U,$$
where $[\pi]$ denotes the $(\pi)^{\infty}$-isotypic
component and where the sum is taken over those $\pi$ such
that
\begin{enumerate}
\item $\pi$ has a $U$-fixed vector,
\item $\pi_{\infty}$ is cohomological for the trivial representation
of $G(\AA^{\infty}_{\QQ})$, and
\item $\pi_q \cong \pi_{\tm,q}$ for all $q$ that split in $E$, are
unramified in $F$, and do not divide the level of $U$. 
\end{enumerate}

This set can be characterized in terms of a near equivalence class
of automorphic representations, in the sense of~\cite{TY}.  We are indebted
to Richard Taylor for the following argument, which establishes this:

\begin{proposition}
For each $\pi$ satisfying conditions (1)-(3), there is a
Galois representation $\rho_{\pi}: \gal(\overline{F}/F) \rightarrow
\GL_n(\overline{\QQ}_l)$ whose restriction to a decomposition group
at almost all $q$ corresponds to $\pi_q$ via the Local Langlands 
Correspondence.  Moreover, if $\rho_{\pi}$ is irreducible 
then $H^{n-1}_{\et}(X_{\overline{F}}, \overline{\QQ}_l)[\pi]$
is isomorphic (up to semisimplification) to some number of copies of $V_{\pi} 
\otimes_{\overline{\QQ}_l} \rho_{\pi^{\prime}}$, where $V_{\pi}$ is
a representation space for $\pi$.
\end{proposition}

\begin{proof}
This is~\cite{TY}, Lemma 3.1, in the case where $\BC(\pi) = (\psi,\Pi)$
with $\JL(\Pi)$ cuspidal, in the notation of~\cite{HT}, VI.2.

In general, $\JL(\Pi)$ is a representation of the form
$$\Pi^{\prime} \boxplus \Pi^{\prime} |\det| \boxplus \dots 
\boxplus \Pi^{\prime} |\det|^{s-1}$$
for some $s$ dividing $n$ and some cuspidal automorphic representation
$\Pi^{\prime}$ of $GL_{\frac{n}{s}}(\AA_F)$.  Since $\JL(\Pi)$ satisfies
$\JL(\Pi)^{\vee} = \JL(\Pi)^c$, we must have 
$(\Pi^{\prime})^{\vee} = (\Pi^{\prime})^c |\det|^{s-1} $.

Fix a character $\psi: \AA_F^{\times}/F^{\times} \rightarrow \CC^{\times}$
with $\psi\psi^c = |\cdot |$.  We can take $\psi$ to be the identity at each
%comment: check this!
archimedean place of $F$.  Then $(\psi \circ \det)\Pi^{\prime}$ 
is a cuspidal automorphic representation of $\GL_{\frac{n}{s}}(\AA_F)$ whose
dual equals its complex conjugate.  Moreover, the representation 
$\JL(\Pi)_{\infty}$
has the same infinitesimal character as an algebraic representation; it
follows that $\Pi^{\prime}_{\infty}$ (and hence 
$((\psi \circ \det)\Pi^{\prime})_{\infty}$) 
do as well.
By~\cite{HT}, VII.1.9 we can associate a Galois representation
$\rho_{(\psi \circ \det)\Pi^{\prime}}$ to 
$(\psi \circ \det)\Pi^{\prime}$.  Taking
direct sums of suitable twists of this representation yields the
desired representation for $\pi$.
\end{proof}

Let $\pi$ and $\pi^{\prime}$ be two representations satisfying the
above conditions.  Then
for $q$ not dividing the level of $U$, unframified in $F$, and split in $E$,
the representations $\rho_{\pi}$ and $\rho_{\pi^{\prime}}$ agree
when restricted to a decomposition group at $q$.  A straightforward
{\v C}ebotarev argument then shows that the two representations are
equal.  It follows that $\pi$ and $\pi^{\prime}$ agree ``almost everywhere'';
i.e. they are nearly equivalent in the sense of~\cite{TY}.  
Thus $\tm$ determines a near equivalence class of automorphic representations.
Moreover, the representation $\rho_{\pi}$ actually depends only on
the near equivalence class of $\pi$, or equivalently, the ideal $\tm$.  
We henceforth denote this representation by $\rho_{\tm}$.

If $m$ is a maximal ideal lying over $\tm$, we also have a reduced
representation 
$\overline{\rho}_m: \gal(\overline{F}/F) \rightarrow GL_n(\TT_U/m)$
via `reduction mod $m$'.  It is well-defined up to semisimplification.

If $\tm^{\prime}$ is another minimal prime of $\TT_U$ contained in $m$,
then $\rho_{\tm^{\prime}}$ will also have the reduced representation
$\overline{\rho}_m$.  Note that if $\overline{\rho}_m$ is absolutely 
irreducible, then $\rho_{\tm^{\prime}}$ will be irreducible as well, and so 
must appear in the cohomology of $X_U$.

Assume $\overline{\rho}_m$ is absolutely irreducible, and
for each $\tm^{\prime}$ contained in $m$, let $V_{\tm^{\prime}}$
be a representation space for $\rho_{\tm^{\prime}}$.  The product
over all $\tm^{\prime}$ of the $V_{\tm^{\prime}}$ is then an
$n$-dimensional representation of $\gal(\overline{F}/F)$ with
coefficients in
$$\prod_{\tm^{\prime}} (\TT_U)_{\tm^{\prime}} \cong (\TT_U)_m \otimes_{\ZZ_l}
\QQ_l.$$
Moreover, when considered as a $(\TT_U)_m \otimes_{\ZZ_l} \QQ_l$-
representation, we have a set of primes of $F$ of density 1 such that
the characteristic polynomial of Frobenius at such primes has
coefficients in $(\TT_U)_m$.  It follows that this representation is
defined over $(\TT_U)_m$; i.e., that we have a representation
$\rho_m: \gal(\overline{F}/F) \rightarrow \GL_n((\TT_U)_m)$,
such that 
$$\rho_m \otimes_{(\TT_U)_m} (\TT_U)_{\tm^{\prime}} = \rho_{\tm^{\prime}}$$
for all $\tm^{\prime}$ contained in $m$.

The upshot is the following, which can be viewed as a result on ``congruences''
for automorphic forms on unitary Shimura varieties:

\begin{corollary} \label{cor:congruences}
Let $U^{\prime}$ be a subgroup containing $U$, and suppose that $m$
descends to a maximal ideal of $\TT_{U^{\prime}}$ via the natural
surjection $\TT_U \rightarrow \TT_{U^{\prime}}$.  Suppose further that
$\overline{\rho}_m$ is absolutely irreducible.  Then the representation
$\overline{\rho}_m$ ``arises from level $U^{\prime}$'', in the sense 
that there is
a miminal prime $\tm^{\prime}$ of $\TT_{U^{\prime}}$ such that
$\rho_m$ is the reduced representation associated to $\rho_{\tm^{\prime}}$.
\end{corollary}

\begin{proof} Since $\TT_{U^{\prime}}$ has no $l$-torsion, there
is a minimal prime $\tm^{\prime}$ of $\TT_{U^{\prime}}$ above $m$.
Then $\rho_{\tm^{\prime}}$ is obtained from $\rho_m$ by tensoring
with $(\TT_U)_{\tm^{\prime}}$.  Since every Jordan-Holder constituent
of $\rho_m$ is isomorphic to $\overline{\rho}_m$, the result follows
immediately.
\end{proof}

%%%%%%%%%%%%%%%%%%%%%%%%%%%%%%%%%%%%%%%%%%%%%%%%%%%%%%%%%%%%%%%%%%%%%%%%%
\section{The weight spectral sequence} \label{sec:weight}

We turn to the study of the representations $\rho_{\tm}$ defined
in the previous section at a decomposition group over $\pp$, in
the case where $U_{\pp}$ is one of the
parahoric subgroups $\Gamma_{i_0,\dots,i_r}$ considered in 
section~\ref{sec:basic}, and the residue characteristic $l$ of $\tm$ is
prime to $p$.  The key is the weight spectral sequence
of Rapoport-Zink.

\begin{theorem}[Rapoport-Zink] \label{theorem:weight}~\cite{RZ2}~\cite{Saito} 
Let $X$ be a proper scheme over $\OO_{F,\pp}$ with strictly semistable
reduction at $\pp$; that is, a regular scheme with smooth general fiber
whose special fiber $Y$ is a divisor with normal crossings.  For each $i$,
let $Y^{(i)}$ denote the disjoint union of all $i+1$-fold intersections of
distinct irreducible components of $Y$.  Then there is a spectral sequence
$$E_1^{p,q} = \oplus_{i\geq max(0,-p)} 
H^{q-2i}_{\et}(Y^{(p+2i)}_{\overline{\FF_p}}, \ZZ_l(-i)) \rightarrow 
H^{p+q}_{\et}(X_{\overline{F_{\pp}}}, \ZZ_l).$$
Modulo torsion, this spectral sequence degenerates at $E_2$.  Moreover
if $N: E_1^{p,q} \rightarrow E_1^{p+2,q-2}(-1)$ denotes the endomorphism 
of this spectral sequence that is the identity on all summands that
appear in both $E_1^{p,q}$ and $E_1^{p+2,q-2}(-1)$ and zero elsewhere,
then $N$ abuts to the monodromy operator on 
$H^{p+q}_{\et}(X_{\overline{F_{\pp}}}, \ZZ_l)$.
\end{theorem}

The local models constructed in section~\ref{sec:reduction} imply
that when $U_{\pp} = \Gamma_{i_0,\dots,i_r}$, the scheme
$(X_U)_{\OO_{F,\pp}}$ has strictly semistable reduction at $\pp$.
Moreover, the nonzero terms of the weight spectral sequence occur
in the range $-r \leq p \leq r$.  We thus have:

\begin{theorem} \label{thm:monodromy} The monodromy operator $N$
satisfies $N^{r+1} = 0$ on $H^{n-1}_{\et}((X_U)_{\overline{F}_{\pp}},\QQ_l)$.
In particular, if $\tm$ is a minimal prime of $\tU$ satisfying the
hyopotheses of section~\ref{sec:representations}, then $N^{r+1}$ is
zero on a representation space for $\rho_{\tm}$.
\end{theorem}

This has consequences for when a mod $m$ Galois representation arising
from $X_U$ can arise from a group with lower level at $\pp$; i.e.
can arise in the cohomology of $(X_U^{\prime})_{\overline{F}_{\pp}}$
for some $U^{\prime}$ containing $U$ and isomorphic to $U$ away from $\pp$, 
but with
$U^{\prime}_{\pp} = \Gamma_{i_0^{\prime}, \dots, i^{\prime}_{r^{\prime}}}$
for some $r^{\prime} < r$.  In particular, a representation 
$\overline{\rho}_m$
cannot arise from such a $U^{\prime}$ unless $N^r$ is zero on a 
representation space for $\overline{\rho}_m$.  A natural question to ask 
is, ``to what extent is the converse true?"

We will give several results and conjectures along these lines, but it
is first necessary to introduce some additional terminology.  Let
$\TT^{\univ}$ be the ``universal'' algebra generated by the Hecke operators
at primes $q$ that split in $E$ and do not divide the level of $U$.
There is a natural surjection $\TT^{\univ} \rightarrow
\TT_U$ which allows us to consider a maximal ideal $m$ of $\TT^U$ as
a maximal ideal of $\TT^{\univ}$.

\begin{definition}
A maximal ideal $m$ of $\TT_U$ is {\em pseudo-Eisenstein} if either:
\begin{itemize}
\item $H^{n-1}_{\et}((X_U)_{\overline{F}_{\pp}}, \ZZ_l)^{\tor}_m$ is nonzero, or
\item $H^i_{\et}((X_U)_{\overline{F}_{\pp}}, \ZZ_l)_m$ is nonzero for some $i$
not equal to $n-1$.
\end{itemize}
\end{definition}

Note that for modular curves, the first of these conditions is vacuous,
while the second condition implies that pseudo-Eisenstein $m$ are actually
Eisenstein, as the Hecke action on $H^0$ and $H^2$ of a modular curve
factors through the Eisenstein ideal.  Pseudo-Eisensteinness thus serves
as an analogue of the Eisenstein condition for unitary Shimura varieties.

The notion of pseudo-Eisensteinness depends on $U$, but note that
if $m$ is pseudo-Eisenstein for $U$, it is automatically pseudo-Eisenstein
for any $U^{\prime}$ contained in $U$.

\begin{conjecture} \label{conj:nilp}
Suppose that $m$ is not pseudo-Eisenstein, that
the residue characteristic $l$ of $m$ is different from $p$, and
that $N^r$ vanishes on a representation space for $\overline{\rho}_m$.  Then
there is a $U^{\prime}$ containing $U$ and isomorphic to $U$ away
from $\pp$, but with 
$U^{\prime}_{\pp} = \Gamma_{i_0^{\prime}, \dots, i^{\prime}_{r^{\prime}}}$
for some $r^{\prime} < r$, such that
$\overline{\rho}_m$ arises from level $U^{\prime}$.
\end{conjecture}

\begin{remark} \rm The pseudo-Eisenstein condition is a bit cumbersome
to verify in practice.  One might hope for a statement along the lines
of `$\overline{\rho}_m$ absolutely irreducible implies $m$ is not 
pseudo-Eisenstein'
but unfortunately this seems quite difficult to prove.  Following
work of Mokrane and Tilouine~\cite{MT}, we give a stronger condition on
$\overline{\rho}_m$ that guarantees $m$ is not pseudo-Eisenstein in the 
Appendix.
\end{remark}

Now
let $n=3$ and fix a compact open $U$ as in section~\ref{sec:basic}, and
assume that $U_{\pp}$ is maximal compact, that is, equal to
$\Gamma_3$.  We denote by $U_{1,2}$ (resp. $U_{2,1}$ and $U_{1,1,1}$)
the subgroups that are equal to $U$ away from $\pp$ and that satisfy 
$(U_{1,2})_{\pp} = \Gamma_{1,2}$, etc.

\begin{theorem} \label{thm:main1}
Let $m$ be a maximal ideal of $\TT_{U_{1,2}}$ (or $\TT_{U_{2,1}}$) that
is not pseudo-Eisenstein and such that $\overline{\rho}_m$ is absolutely 
irreducible.
Suppose that $\overline{\rho}_m$ is unramified at $\pp$;
i.e., that $N = 0$ on a representation space for $\overline{\rho}_m$, and that
$\overline{\rho}_m(\Frob_{\pp})$ has distinct eigenvalues.  Then 
$\overline{\rho}_m$ arises from level $U$. 
\end{theorem}

\begin{theorem} \label{thm:main2}
Let $m$ be a maximal ideal of $\TT_{U_{1,1,1}}$ that is not 
pseudo-Eisenstein and
such that $\overline{\rho}_m$ is absolutely irreducible.  Suppose
that $N^2 = 0$ on a representation space for $\overline{\rho}_m$, and
that $l$ does not divide $p^{2d_{\pp}} - 1$.  Then $\overline{\rho}_m$ arises
from level $U_{1,2}$ (or, equivalently, from level $U_{2,1}$.)
\end{theorem}

%%%%%%%%%%%%%%%%%%%%%%%%%%%%%%%%%%%%%%%%%%%%%%%%%%%%%%%%%%%%%%%%%%%%%%%%%
\section{Proof of the main theorems} 
\label{sec:degeneration}

We will now use a more detailed analysis of the weight spectral
sequence to prove Theorems~\ref{thm:main1} and~\ref{thm:main2}.

First consider Theorem~\ref{thm:main1}.  Let $m$ be a maximal
ideal of $\TT_{U_{1,2}}$, and assume that $m$ is not pseudo-Eisenstein,
$\overline{\rho}_m$ is absolutely irreducible, and $\overline{\rho}_m$ is 
unramified at $\pp$.

We wish to show that if $\overline{\rho}_m(\Frob_{\pp})$ has distinct 
eigenvalues, then $\overline{\rho}_m$ arises as a Jordan-holder consituent
of $H^{n-1}_{\et}((X_U)_{\overline{F}}, \ZZ_l)$.  In fact, it suffices to show
that $H^{n-1}_{\et}((X_U)_{\overline{F}}, \ZZ_l)_m$ is nonzero.  For this
space is non-torsion (since $m$ is not pseudo-Eisenstein), and therefore
$m$ descends to a maximal ideal of $\TT_U$.  By Corollary~\ref{cor:congruences},
$\overline{\rho}_m$ then arises from level $U$.

Assume that $H^{n-1}_{\et}((X_U)_{\overline{F}}, \ZZ_l)_m = 0$.  We will
show that $\overline{\rho}_m(\Frob_{\pp})$ has a multiple eigenvalue.  
Consider the weight spectral sequence for $X_{U_{1,2}}$, localized
at $m$.  It has the following form:

$$
\begin{array}{ccccc}
H^2((X_U)^{(1)}_{\overline{\FF}_p}, \ZZ_l(-1))_m & \rightarrow &
H^4(((X_{U_{1,2}})_0)_{\overline{\FF}_p}, \ZZ_l)_m \oplus
H^4(((X_{U_{1,2}})_1)_{\overline{\FF}_p}, \ZZ_l)_m\\
H^1((X_U)^{(1)}_{\overline{\FF}_p}, \ZZ_l(-1))_m & \rightarrow &
H^3(((X_{U_{1,2}})_0)_{\overline{\FF}_p}, \ZZ_l)_m \oplus
H^3(((X_{U_{1,2}})_1)_{\overline{\FF}_p}, \ZZ_l)_m\\
H^0((X_U)^{(1)}_{\overline{\FF}_p}, \ZZ_l(-1))_m & \rightarrow &
H^2(((X_{U_{1,2}})_0)_{\overline{\FF}_p}, \ZZ_l)_m \oplus
H^2(((X_{U_{1,2}})_1)_{\overline{\FF}_p}, \ZZ_l)_m &
\rightarrow & H^2((X_U)^{(1)}_{\overline{\FF}_p}, \ZZ_l)_m\\
& &
H^1(((X_{U_{1,2}})_0)_{\overline{\FF}_p}, \ZZ_l)_m \oplus
H^1(((X_{U_{1,2}})_1)_{\overline{\FF}_p}, \ZZ_l)_m &
\rightarrow & H^1((X_U)^{(1)}_{\overline{\FF}_p}, \ZZ_l)_m\\
& &
H^0(((X_{U_{1,2}})_0)_{\overline{\FF}_p}, \ZZ_l)_m \oplus
H^0(((X_{U_{1,2}})_1)_{\overline{\FF}_p}, \ZZ_l)_m &
\rightarrow & H^0((X_U)^{(1)}_{\overline{\FF}_p}, \ZZ_l)_m\\
\end{array}
$$

The three left-hand maps are Gysin maps; the three right-hand maps
are restrictions.  Note that since $(X_{U_{1,2}})_0$ is isomorphic
to $(X_U)_{\OO_F/\pp}$, and $X_U$ has good reduction at $\pp$,
$H^i(((X_{U_{1,2}})_0)_{\overline{\FF}_p}, \ZZ_l)_m$ is isomorphic
to $H^i((X_U)_{\overline{F}_{\pp}}, \ZZ_l)_m$ and is therefore zero
because of our assumptions on $m$.

\begin{lemma} \label{lemma:gysin}
The Gysin maps 
$$H^i((X_U)^{(1)}_{\overline{\FF}_p}, \ZZ_l(-1))_m \rightarrow
H^{i+2}(((X_{U_{1,2}})_1)_{\overline{\FF}_p}, \ZZ_l)_m$$
all vanish.
\end{lemma}
\begin{proof}
To shorten notation, let $\tX$, $X$, and $Y$ denote the
schemes $((X_{U_{1,2}})_1)_{\overline{\FF}_p}$,
$(X_U)_{\overline{\FF}_p}$, and $(X_U)^{(1)}_{\overline{\FF}_p},$
respectively.
We will show that the composition of the Gysin map:
$$H^i(Y, \ZZ_l(-1)) \rightarrow
H^{i+2}(X, \ZZ_l)$$
with the natural map:
$$H^{i+2}(X, \ZZ_l) \rightarrow H^{i+2}(\tX, \ZZ_l)$$
is equal to $p^{d_{\pp}}$ times the Gysin map:
$$H^i(Y, \ZZ_l(-1)) \rightarrow
H^{i+2}(\tX, \ZZ_l).$$
The result follows immediately since the cohomology of
$X$ is isomorphic to the cohomology of
$(X_U)_{\overline{F}}$ (as the latter has good reduction at $\pp$)
and hence by assumption the cohomology of $X$ vanishes after 
localizing at $m$.

Let $\ti$ denote the inclusion of $Y$ in
$\tX$ (as the intersection of $(X_{U_{1,2}})_0$ with
$(X_{U_{1,2}})_1$), and let $i$ denote the inclusion of
$Y$ in $X$.  Let $j$ be the natural map
$\tX \rightarrow X$.  Then the Gysin maps
under consideration are induced by the isomorphisms
$R^2\ti^!(\ZZ_l)_{\tX} \cong \ZZ_l(-1)_Y$ and 
$R^2i^!(\ZZ_l)_X \cong \ZZ_l(-1)_Y$, where $(\ZZ_l)_{\tX}$, $(\ZZ_l)_X$,
and $(\ZZ_l(-1))_Y$ are constant sheaves on $\tX$, $X$, and $Y$, 
respectively.
The assertion of the previous paragraph thus amounts to the assertion that
in the commutative diagram:
$$
\begin{array}{ccc}
R^2i^!(\ZZ_l)_X & \rightarrow & R^2\ti^!j^*(\ZZ_l)_X\\
\downarrow & & \downarrow\\
\ZZ_l(-1)_Y & \rightarrow & \ZZ_l(-1)_Y
\end{array}
$$
the bottom arrow is multiplication by $p^{d_{\pi}}$.  This assertion
is {\'e}tale local in nature, and can be verified generically on $Y$.

Let $x$ be a point of $Y$ outside $(X_U)^{(0)}.$  Then in an
{\'e}tale neighborhood of $\ti(x)$, the map $\tX \rightarrow X$
looks like the composition of a purely inseparable morphism of
degree $p^{d_{\pi}}$ with an {\'e}tale map $\tX \rightarrow X$.
Moreover, the map $\tX \rightarrow X$ induces an isomorphism
of $\ti(Y)$ with $i(Y)$.  Since $\tX$, $X$, and $Y$ are all smooth,
this means that in an {\'e}tale neighborhood of $\ti(x)$, the map
$\tX \rightarrow X$ looks like the map 
$\AA^2_{\overline{\FF}_p} \rightarrow \AA^2_{\overline{\FF}_p}$ that
raises the second coordinate to the $p^{d_{\pp}}$th power.

It thus suffices to show that the map this induces 
$$H^2_{\AA^1_{\overline{\FF}_p}}(\AA^2_{\overline{\FF}_p}, \ZZ_l)
\rightarrow
H^2_{\AA^1_{\overline{\FF}_p}}(\AA^2_{\overline{\FF}_p}, \ZZ_l)$$
is multiplication by $p^{d_{\pi}}$.  (Here $\AA^1$ is considered
as the subscheme of $\AA^2$ consisting of points with second coordinate
zero.)  But the long exact sequence for cohomology with support identifies
this with the map
$$H^1((\AA^2 \setminus \AA^1)_{\overline{\FF}_p}, \ZZ_l) \rightarrow
H^1((\AA^2 \setminus \AA^1)_{\overline{\FF}_p}, \ZZ_l)$$
induced by raising the second coordinate to the $p^{d_{\pp}}$th power;
this map is clearly multiplication by $p^{d_{\pp}}$.
\end{proof}

It follows that the maps in the left-hand column of the weight
spectral sequence for $X_{U_{1,2}}$ vanish after localization at $m$.
By Poincar{\'e} duality it follows that the right-hand maps all have
torsion image, but since every cohomology group appearing in the
right-hand column is torsion-free, the right-hand maps vanish as well.
Thus the $E_2$ terms of the spectral sequence are the same as the
$E_1$ terms.  Moreover, the weight spectral sequence degenerates at $E_2$
modulo torsion, and the only remaining nonzero boundary maps map
into cohomology groups that are torsion-free.  Thus, after localization
at $m$, the weight spectral sequence is already degenerate at $E_1$.

On the other hand, this spectral sequence abuts to
$H^{p+q}((X_{U_{1,2}})_{\overline{F}_{\pp}}, \ZZ_l)_m$, which is
zero except in degree 2.  It follows that all terms of the above
spectral sequence vanish except in cohomological degree $2$.  We obtain
a filtration of $H^2((X_{U_{1,2}})_{\overline{F}_{\pp}}, \ZZ_l)_m$
whose successive quotients are
$E_{\infty}^{1,1} = H^1((X_U)^{(1)}_{\overline{\FF}_p}, \ZZ_l)_m$,
$E_{\infty}^{0,2} = H^2(((X_{U_{1,2}})_1)_{\overline{\FF}_p}, \ZZ_l)_m$, and
$E_{\infty}^{-1,3} = H^1((X_U)^{(1)}_{\overline{\FF}_p}, \ZZ_l(-1))_m$.  
Up to a Tate twist,
the monodromy operator induces an isomorphism of $E_{\infty}^{-1,3}$ 
with $E_{\infty}^{1,1}$, and is zero on the remaining quotients.

\begin{remark} \label{remark:support} \rm Note that the above argument
in fact shows that every $m$ in the support of
$H^i((X_U)^{(1)}_{\overline{\FF}_p}, \ZZ_l)$ is either pseudo-Eisenstein,
or in the support of $H^2((X_{U_{1,2}})_{\overline{F}_{\pp}}, \ZZ_l)$.
Similarly, every $m$ in the support of
$H^i(((X_{U_{1,2}})_1)_{\overline{\FF}_p}, \ZZ_l)$ is either
pseudo-Eisenstein or in the support of
$H^2((X_{U_{1,2}})_{\overline{\FF}_p}, \ZZ_l)$.
\end{remark}

\begin{remark} \label{remark:compgrp} \rm 
A key step in the classical proof of Mazur's principle
for modular curves is showing that the maximal ideal $m$ under consideration
is not in the support of the component group of the Jacobian of the modular
curve.  After tensoring with $\ZZ_l$, we can think of this component group 
as the cokernel of the monodromy operator $N$ (considered as a map from the
highest weight quotient of the middle cohomology of the modular curve to
the lowest weight submodule of this cohomology.)  The claim that the 
component group vanishes when localized at $m$ is thus equivalent
to the claim that this $N$ is surjective.
The fact that in our setting the monodromy operator induces
an isomorphism of $E_{\infty}^{-1,3}$ with $E_{\infty}^{1,1}$ 
after localizing at $m$ is therefore an analogue of this statement.

Note that even in the classical case, this ``Eisensteinness of
the component group'' need not hold if the modular curve under consideration
is not a fine moduli space (see for instance~\cite{Ri2}.)  Thus our requirement
that $U$ be ``sufficiently small'' is more than just a technical convenience.
\end{remark}

Now consider the representation $\rho_m$.  Since $\overline{\rho}_m$
does not arise from level $U$, each prime $\tm$ of $\TT_{U_{2,1}}$
contained in $m$ corresponds to a near equivalence class of automorphic 
representations $\pi_{\tm}$ that have no $U$-fixed vectors.  It follows
by the compatibility of local and global Langlands that all such
$\rho_{\tm}$ are ramified at $p$.  If $V_{\tm}$ is a representation
space for $\rho_{\tm}$, then $V_{\tm}$ is a three-dimensional
$(\TT_{U_{2,1}})_{\tm}$-vector space, equipped with a nonzero operator
$N$ such that $N^2 = 0$.  The successive quotients of the monodromy filtration 
on $V_{\tm}$ are therefore the one dimensional spaces $\im N$, $\ker N/\im N$,
and $N/\ker N$.  

It therefore follows that if $V_m$ is a representation space for $\rho_m$,
then $\im N$, $\ker N/\im N$, and $N/\ker N$ are free $(\TT_{U_{1,2}})_m$-
modules of rank one, and therefore $\rho_m(\Frob_{\pp})$ acts on them
via scalars $\alpha_3$, $\alpha_2$, and $\alpha_1$, respectively.
(The subscripts correspond to the weights of the Frobenius action.)
Therefore $\Frob_{\pp}$ must act on the successive quotients 
$E_{\infty}^{1,1}$ $E_{\infty}^{0,2}$ and $E_{\infty}^{-1,3}$  
of the monodromy filtration on 
$H^2((X_{U_{1,2}})_{\overline{F}_{\pp}}, \ZZ_l)_m$ 
by the same scalars $\alpha_3,\alpha_2,\alpha_1$.
As we have seen, these quotients are
$H^1((X_U)^{(1)}_{\overline{\FF}_p}, \ZZ_l)_m$,
$H^2(((X_{U_{1,2}})_1)_{\overline{\FF}_p}, \ZZ_l)_m$, and
$H^1((X_U)^{(1)}_{\overline{\FF}_p}, \ZZ_l(-1))_m$, respectively.

Now since $\overline{\rho}_m$ is unramified at $\pp$, $N$ is zero
on $H^2((X_{U_{1,2}})_{\overline{F}_{\pp}}, \ZZ_l)_m$ modulo $m$.  It
follows that $E_{\infty}^{1,1}$ is contained in
$mH^2((X_{U_{1,2}})_{\overline{F}_{\pp}}, \ZZ_l)_m$.  In particular
the quotient $M$ of $H^2((X_{U_{1,2}})_{\overline{F}_{\pp}}, \ZZ_l)_m$
by $E_{\infty}^{1,1}$ surjects onto 
$H^2((X_{U_{1,2}})_{\overline{F}_{\pp}}, \ZZ_l)_m \otimes_{\TT_{U_{1,2}}}
\TT_{U_{1,2}}/m.$

But the latter is a direct sum of copies of a representation space
$V_{\overline{\rho}_m}$ for $\overline{\rho}_m$.  In particular
$\alpha_3$, $\alpha_2$, and $\alpha_1$ all appear (modulo $m$) as 
eigenvalues of Frobenius on this space.  On the other hand $M$ has
a two step filtration such that $\Frob_{\pp}$ acts via $\alpha_2$ and
$\alpha_3$ on the two successive quotients.  Thus $\alpha_1$
must be congruent to either $\alpha_2$ or $\alpha_3$ modulo $m$.
But $\alpha_3$, $\alpha_2$, and $\alpha_1$ reduce modulo $m$ to 
the eigenvalues of
$\overline{\rho}_m(\Frob_{\pp})$.
and so $\overline{\rho}_m$ has a multiple eigenvalue.  We have
thus proved Theorem~\ref{thm:main1}.
 
We now turn to Theorem~\ref{thm:main2}.  Take $m$ to be a maximal
ideal of $\TT_{U_{1,1,1}}$, and assume that $m$ is not pseudo-Eisenstein,
$\overline{\rho}_m$ is absolutely irreducible, and $N^2 = 0$ on
a representation space for $\overline{\rho}_m$.

As before, it suffices to show that either 
$H^2((X_{U_{1,2}})_{\overline{F}_{\pp}}, \ZZ_l)_m$ is nonzero
or that $l$ divides $p^{2d_{\pp}} - 1$.  We assume that 
$H^2((X_{U_{1,2}})_{\overline{F}_{\pp}}, \ZZ_l)_m$ vanishes,
and will show that this implies that $l$ divides $p^{2d_{\pp}} - 1$.

This assumption, together with Remark~\ref{remark:support}, shows that
after localizing at $m$, the only nonzero terms in the weight
spectral sequence for $X_{U_{1,1,1}}$ are
$$E_{\infty}^{-2,4} = H^0((X_U)^{(2)}_{\overline{\FF}_p}, \ZZ_l(-2))_m,$$
$$E_{\infty}^{0,2} = H^0((X_U)^{(2)}_{\overline{\FF}_p}, \ZZ_l(-1))_m,$$
$$E_{\infty}^{2,0} = H^0((X_U)^{(2)}_{\overline{\FF}_p}, \ZZ_l)_m.$$
The operator $N^2$ identifies $E_{\infty}^{-2,4}$ with $E_{\infty}^{2,0}(-2)$.

On the other hand, if $\tm$ is a prime of $\TT_{U_{1,1,1}}$ containing
$m$, then by our assumtpion $\tm$ corresponds to a near equivalence
class of representations with a $U_{1,1,1}$-fixed vector but no
$U_{1,2}$ or $U_{2,1}$-fixed vectors.  Taylor-Yoshida show~\cite{TY} that
for such representations, the mondromy operator on a representation
space for $\rho_{tm}$ satisfies $N^2 \neq 0$.  In particular $N$
has rank two, and the successive quotients on the monodromy filtration
on a representation space for $\rho_{\tm}$ are $\im N^2$, $\im N/\im N^2$,
and $N/\im N$.  

It follows as in the previous case that $\Frob_{\pp}$ acts
on $E_{\infty}^{-2,4}$, $E_{\infty}^{0,2}$, and $E_{\infty}^{2,0}$
by scalars $\alpha_4$, $\alpha_2$, and $\alpha_0$ in $(\TT_{U_{1,1,1}})_m$.
Moreover, since these spaces are simply Tate twists of one another, we
have $$\alpha_4 = p^{d_{\pi}}\alpha_2 = p^{2d_{\pi}}\alpha_0.$$

Now the argument proceeds exactly as in the previous case.  We
have that 
$$E_{\infty}^{2,0} \subset 
mH^2((X_{U_{1,1,1}})_{\overline{F}_{\pp}}, \ZZ_l)_m,$$
so that
if $M$ is the quotient $H^2((X_{U_{1,1,1}})_{\overline{F}_{\pp}}, \ZZ_l)_m/
E_{\infty}^{2,0}$, then $M/mM$ contains a copy of $\overline{\rho}_m$.
Thus $\alpha_0$ must be congruent to either $\alpha_2$ or $\alpha_4$
modulo $m$.  This can only happen if $p^{2d_{\pi}}$ is congruent to
one modulo the residue characteristic $l$ of $m$, as was to be proven.
%%%%%%%%%%%%%%%%%%%%%%%%%%%%%%%%%%%%%%%%%%%%%%%%%%%%%%%%%%%%%%%%%%%%%%%%%
\section{Appendix: A criterion for non-pseudo-Eisensteinness}

Since the non-pseudo-Eisenstein condition of section~\ref{sec:weight} is
cumbersome, it would be useful to give conditions that imply it, but are
more concrete and easier to verify in practice.  Here, we give one
such criterion, which is a $\GL_n$ analogue of Theorem 1 of~\cite{MT}
(which deals with Siegel modular varieties).  The argument here is
a straightforward adaptation of that given in~\cite{MT}; we reproduce it
here mainly for self-containedness, and also to point out the (slight)
changes necessary for the $\GL_n$ case.  Note that since we only work
with constant coefficient-sheaves, the complicated technical machinery
of Bernstein-Gelfand-Gelfand complexes employed by Mokrane and
Tilouine is not needed here.

Consider a (fixed) unitary Shimura variety $X_U$, of the form we have
been considering, and let $m$ be a maximal ideal of $\TT_U$, of
residue characteristic $l > 5$.  We impose the following conditions on
the mod $m$ representation $\rho_m$:

\noindent
{\em Ordinarity}- There is a prime $\ll$ of $\OO_F$ above $l$, such
that the restriction of $\overline{\rho}$ to a decomposition group at
$\ll$ has the form:
$$
\begin{pmatrix}
\chi^{n-1} & * & \dots & * \\
0 & \chi^{n-2} & \dots & * \\
  & \vdots & \\
0 & 0 & \dots & 1 \\
\end{pmatrix},
$$
where $\chi$ denotes the mod $l$ cyclotomic character. We require that
the characters $\chi^0, \chi^1, \dots, \chi^{n-1}$ induce
distinct maps $I_{\ll} \rightarrow \FF_l^{\times}$, where $I_{\ll}$
is the inertia group at $\ll$.

\noindent
{\em Large Image-} There exists a split and reductive (but not necessarily
connected) subgroup $H$ of $\GL_n$, and a subfield $k^{\prime}$ of
$\TT_U/m$ such that:
\begin{itemize}
\item $H$ contains a maximal torus $T$ of $\GL_n$, 
\item $H$ contains a subgroup $W'$ of the Weyl group of $T$ 
that acts transitively on the weights of the standard representation
of $\GL_n$,
\item the image of the inertia group $I_{\ll}$ under $\overline{\rho}_m$
is contained in $H^0(k^{\prime})$ (here $H^0$ is the connected component
of the identity in $H$), and
\item the image of $\overline{\rho}_m$ contains the subgroup
$H(k^{\prime})_{\det}$ of $H(k^{\prime})$ consisting of those elements of
$H(k^{\prime})$ whose determinant is in the image of $\det \overline{\rho}_m$.
\end{itemize}

\begin{remark} \rm The ``Large Image'' condition appearing here (and its
analogue in~\cite{MT}) really says that the image of $\overline{\rho}_m$
is ``large relative to the image of $I_{\ll}$''.  In particular, when
$n=3$, if the image of $I_{\ll}$ s not contained in the $k^{\prime}$-points
of any maximal torus
$T$ of $\GL_n$, then the only possible choice for $H$ is $H = \GL_3$.
On the other hand, if the image of $I_{\ll}$ is contained in $T(k^{\prime})$
for some $T$ and $k^{\prime}$, one can take $H$ to be the semidirect 
product of $T$ with any suitable $W^{\prime}$.
\end{remark}

Under these assumptions, we find:
\begin{theorem} \label{thm:pseudo}
The maximal ideal $m$ is not pseudo-Eisenstein;
in particular:
\begin{itemize} 
\item $H^{n-1}_{\et}((X_U)_{\overline{F}}, \ZZ_l)^{\tor}_m = 0$, and
\item $H^i_{\et}((X_U)_{\overline{F}}, \ZZ_l)_m = 0$ for $i \neq n-1$.
\end{itemize}
\end{theorem}

To prove this, we follow~\cite{MT}, sections 7.1 and 7.2.  We begin
by studying representations of subgroups $H$ of the sort arising
in the ``Large Image'' condition.  Fix an $H$, $T$, $W^{\prime}$
and $k^{\prime}$ as described in that condition.

We follow the notation in~\cite{MT} quite closely.  We let $X$ denote
the character group of $T$, let $\Psi_{H^0}$ denote the roots of
$H^0$, and $\Psi$ denote the roots of $\GL_n$.  Denote by
$\Psi_{H^0}^{\perp}$ the subspace of those $\lambda$ in $X$ such that
$\<\lambda,\beta^{\vee}\> = 0$ for all $\beta$ in $\Psi_{H^0}$,
where $\beta^{\vee}$ is the coroot associated to $\beta$.
Choose a set of positive roots $\Psi^+$ for $\GL_n$; the corresponding
Borel subgroup determines a Borel subgroup for $H^0$ as well, and hence
a set $\Psi_{H^0}^+$ of positive roots for $H^0$.  Let 
$\Delta_{H^0}$ denote the set of positive simple roots for $H^0$.

Let $T(k^{\prime})_{\det}$ be the subgroup of $T(k^{\prime})$ consisting
of those elements whose determinant is in the image of
$\det \overline{\rho}_m$.  We have a map
$$X \rightarrow \Hom(T(k^{\prime})_{\det}, (k^{\prime})^{\times})$$
whose kernel is the sublattice
$$(q^{\prime} - 1)X + e \cdot \ZZ \cdot \det,$$
where $q^{\prime}$ is the order of $k^{\prime}$ and $e$ is the index
of the image of $\det \overline{\rho}_m$ in $k^{\prime}$.  
Let $\tPsi_{H^0}^{\perp}$ denote the space
$(q-1)\Psi_{H^0}^{\perp} + e \cdot \ZZ \cdot \det$; that is, the
intersection of $\Psi_{H^0}^{\perp}$ with this kernel.

Let $X_{H,q^{\prime}}$ denote the set
$$\{(v,a) \in X/\tPsi_{H^0}^{\perp} : 
0 \leq \<v, \beta^{\vee}\> \leq q^{\prime} - 1 \forall \beta \in 
\Delta_{H^0}\}.$$ 
By Steinberg's theorem, (c.f.~\cite{MT}, p. 52),
the irreducible representations of $H^0(k^{\prime})_{\det}$, as an abstract
group, are classified by elements of $X_{H,q^{\prime}}$.  In particular,
if $V$ is the irreducible algebraic representation of $H^0$ with highest
weight $\mu$, and $\mu$ is an element of $X_{H,q^{\prime}}$ (when
considered modulo $\tPsi_{H^0}^{\perp}$) 
then $V_{k^{\prime}}$ is an irreducible representation
of $H^0(k^{\prime})_{\det}$ as an abstract group, and this representation
depends only on the class of $\mu$ modulo $\tPsi_{H^0}^{\perp}$.

\begin{lemma} (c.f.~\cite{MT}, Lemma 12)
Let $\mu$ be a dominant weight of $H^0$, and $V$ the
irreducible algebraic representation of $\mu$ with highest weight
$\mu$.  Suppose that:
\begin{itemize}
\item $V_{k^{\prime}}$ is an irreducible representation of
$H^0({k^{\prime}})_{\det}$,
\item when considered as a map $T(k^{\prime})_{\det} \rightarrow 
(k^{\prime})^{\times}$, $\mu$ is equal to the highest weight of
the standard representation of $G$, and
\item every weight of $V$ coincides with some weight of the standard
representation of $G$ when both are considered as maps
$T(k^{\prime})_{\det} \rightarrow (k^{\prime})^{\times}$.
\end{itemize}
Then $V_{k^{\prime}}$ coincides with the ``standard representation''
of $H^0({k^{\prime}})_{\det}$; i.e., $\mu$ is equal to the highest weight
of the standard representation modulo $\tPsi_{H^0}^{\perp}$.
\end{lemma}
\begin{proof}
Let $w$ denote the highest weight of the standard representation.  By
hypothesis, $\mu-w \in (q^{\prime} - 1)X + e \cdot \ZZ \cdot \det$.  It thus
suffices to show that $\mu - w$ also lies in $\Psi_{H^0}^{\perp}$, as
then it must lie in $\tPsi_{H^0}^{\perp}$ as required.

Observe that for any positive simple root $\alpha$ of $H^0$ we must
have $\<w,\alpha\> \in \{-1,0,1\}$.  Moreover, by our assumptions
we have
$0 \leq \<\mu, \alpha\>  \leq q^{\prime} - 1$ and
$q^{\prime} - 1$ must divide $\<w - \mu, \alpha\>$ for all $\alpha$.

If $\<w,\alpha\> = 1$, these observations immediately imply that
$\<\mu,\alpha\>$ is also equal to $1$.

On the other hand, if $\<w,\alpha\> = 0$, a priori we could have
$\<\mu,\alpha\> = 0$ or $\<\mu,\alpha\> = q^{\prime} - 1$.  In the
latter case, however, the $\alpha$-string of $\mu$ would have
length $q^{\prime} - 1$, and so $\mu - \alpha$ is also a weight of
$V$.  Thus there is a weight $w^{\prime}$ of the standard
representation that coincides with $\mu - \alpha$ when considered as a map
$T(k^{\prime})_{\det} \rightarrow (k^{\prime})^{\times}$.  We then
have $\<w^{\prime}, \alpha\> \in \{-1,0,1\}$, and $q^{\prime} - 1$
divides $\<w^{\prime} - (\mu - \alpha), \alpha\>$.  On the other
hand, $\<\mu - \alpha, \alpha\> = q^{\prime} - 3$, so this is only
possible if $q^{\prime} - 1$ divides $1$, $2$, or $3$.  This cannot
occur since $q^{\prime}$ is a power of $l$.

If $\<w, \alpha\> = -1$, then $\<\mu, \alpha\> = q^{\prime} - 2$.
Thus $\mu - \alpha$ is a weight of $V$.  By the same argument as above,
we find that in this case $q^{\prime} - 1$ must divide $2$, $3$, or $4$,
and none of these can occur.  Thus $\<w,\alpha\>$ is never equal to
$-1$.
\end{proof}

\begin{corollary} Let $V$ be a nonzero $\TT_U/m$-vector space with an
action of $\gal(\overline{F}/F)$-module, and
suppose that for each $\sigma \in \gal(\overline{F}/F)$,
$P_{\sigma}(\sigma)$ annihilates $V$, where $P_{\sigma}$ is the
characteristic polynomial of $\overline{\rho}_m(\sigma)$.  Then
$V$ has a subquotient on which $I_{\ll}$ acts via $\overline{\rho}_m$.
\end{corollary}
\begin{proof}
Let $\Gamma^{\prime}$ denote the preimage under $\overline{\rho}_m$
of $H^0(k^{\prime})_{\det}$ in $\gal(\overline{F}/F)$, and observe
that $\Gamma^{\prime}$ contains $I_{\ll}$ by our ``Large Image'' assumption.
Let $\Gamma^{\prime\prime}$ denote the kernel of the map
$\Gamma^{\prime} \rightarrow H^0(k^{\prime})_{\det}$ induced by
$\overline{\rho}_m$.  Then $\Gamma^{\prime\prime}$ is a nilpotent
$p$-group, so without loss of generality (replacing $V$ with the 
$\Gamma^{\prime\prime}$-invariants of $V$, which are necessarily nontrivial)
we can assume $V$ is fixed by $\Gamma^{\prime\prime}$.  Thus we can
view $V$ as a $H^0(k^{\prime})_{\det}$-module, by identifying
$H^0(k^{\prime})_{\det}$ with $\Gamma^{\prime}/\Gamma^{\prime\prime}$.

Any irreducible $H^0(k^{\prime})_{\det}$-constituent of $V$ arises
(as above) from an irreducible $H^0$-representation $V^{\mu}$ with
some highest weight $\mu$.  On the other hand, by our characteristic
polynomial assumption, the only maps $T(k^{\prime})_{\det} \rightarrow
(k^{\prime})^{\times}$ appearing in the decomposition of
$V$ into irreducible $T(k^{\prime})_{\det}$-modules are the characters
of the standard representation of $\GL_n$.  Thus for every weight of 
$V^{\mu}$ there is a character of the standard representation that
induces the same map $T(k^{\prime})_{\det} \rightarrow 
(k^{\prime})^{\times}.$  Moreover, since $W^{\prime}$ acts on $V$, and
acts transitively on the weights of the standard representation, we
can assume (replacing $\mu$ with $\mu^y$ for a suitable $y \in W^{\prime}$)
that $\mu$ induces the same map $T(k^{\prime})_{\det} \rightarrow
(k^{\prime})^{\times}$ as the {\em highest} weight of the
standard representation of $\GL_n$.  It then follows from the
previous lemma that $V^{\mu}_{k^{\prime}}$ is the ``standard representation''
of $H^0(k^{\prime})_{\det}$, and is an irreducible constituent of $V$.

But this means precisely that $\Gamma^{\prime}$ acts on $V^{\mu}_{k^{\prime}}$
via $\overline{\rho}_m$.  In particular $I_{\ll}$ acts via $\overline{\rho}_m$,
as desired.
\end{proof}

Theorem~\ref{thm:pseudo} now follows easily.  As in~\cite{MT} it
suffices to show that $H^i_{\et}((X_U)_{\overline{F}}, \FF_l)_m$
vanishes for $i < n-1$.  Such a cohomology group is annihilated
by $P_{\sigma}(\sigma)$ for all $\sigma$ in $\gal(\overline{F}/F)$.
The preceding corollary thus shows that, if nonzero, it contains a subquotient
on which $I_{\ll}$ acts via $\overline{\rho}_m$; by ordinarity
it would follow that $H^i_{\et}((X_U)_{\overline{F}}, \FF_l)_m$ had
at least $n$ distinct Hodge-Tate weights.  But the
{\'e}tale-DeRham comparison theorem shows that the Hodge filtration
on $H^i_{\et}((X_U)_{\overline{F}}, \FF_l)$ has length $i+1$, and,
as $i+1 < n$, $H^i_{\et}((X_U)_{\overline{F}}, \FF_l)_m$ must vanish as 
required.

\textsc{Acknowledgements}

The author is very grateful to Richard Taylor for his advice and
encouragement.  The work on this paper was partially supported by the 
National Science Foundation.

\end{document}